\documentclass{article}
\usepackage{amssymb}

\usepackage{graphicx}

\newtheorem{problem}{Problem}[section]
\newtheorem{theo}[problem]{Theorem}
\newtheorem{rem}[problem]{Remark}
\newtheorem{prob}[problem]{Problem}
\newtheorem{defin}[problem]{Definition}
\newtheorem{prop}[problem]{Proposition}
\newtheorem{cor}[problem]{Corollary}
\newtheorem{lema}[problem]{Lemma}
\newtheorem{exam}[problem]{Example}
\newtheorem{observ}[problem]{Observation}

\begin{document}
\date{October  2003}
\title{{\Large Topology and Combinatorics of\\
Partitions of Masses by Hyperplanes}}

\vskip 2cm

\author{{\Large Peter Mani--Levitska\thanks{Supported
by Schweizerischer Nationalfonds zur F\" orderung der
wissenschaftlichen Forschung.}}
\\ {Mathematics Institute, Bern} \\[5mm]
{\Large Sini\v sa Vre\' cica\thanks{Supported by the Ministry for
Science and Technology of Serbia, Grant 1643.}}
\\ {Faculty of Mathematics, Belgrade}\\[5mm]
  {\Large  Rade \v Zivaljevi\' c\thanks{Supported by the Ministry for
Science and Technology of Serbia, Grant 1643.}}
 \\ {Mathematics Institute
 SANU, Belgrade}}
\maketitle

\vskip 2cm
\begin{abstract}
An old problem in combinatorial geometry is to determine when one
or more measurable sets in $\mathbb{R}^d$ admit an {\em
equipartition} by a collection of $k$ hyperplanes, \cite{Gru}
\cite{Hadw}. The problem can be reduced to the question of
(non)existence of a map $f : (S^d)^k\rightarrow S(U)$, equivariant
with respect to the Wayl group $W_k := ({\mathbb{Z}/2})^{\oplus
k}\rtimes S_k$, where $U$ is a representation of $W_k$ and
$S(U)\subset U$ the corresponding unit sphere. In this paper we
develop a general method for computing topological obstructions
for the existence of such equivariant maps. Emphasizing the
combinatorial point of view, we show that the computation of
relevant cohomology/bordism obstruction classes can be in many
cases reduced to the question of enumerating the classes of
immersed curves in $\mathbb{R}^2$ with a prescribed type and
number of intersections with the coordinate axes. It turns out
that the last problem is closely related to some ``cyclic word
enumeration problems'' from enumerative combinatorics which are
usually approached via Poly\' a enumeration or M\" obius inversion
technique. Among the new results is the well known open case of
$5$ measures and $2$ hyperplanes in $\mathbb{R}^8$, \cite{Ram}.
The obstruction in this case is identified as the element $2
X_{ab}\in H_1(\mathbb{D}_8; {\cal Z})\cong \mathbb{Z}/4$, where
$X_{ab}$ is a generator, which explains  why this result cannot be
obtained by the parity count formulas from \cite{Ram} or the
methods based on either Stiefel-Whitney classes or ideal valued
cohomological index theory  \cite{FH2}.
\end{abstract}

\vfill\newpage

\tableofcontents

\vfill\newpage

\section{Introduction}
\label{Intro}

\subsection{The Equipartition Problem}
\label{Pocetak}

\begin{defin}
\label{equi} Suppose that $${\cal M} = \{\mu_1, \mu_2, \ldots ,
\mu_j\}$$ is a collection of continuous mass
distributions/measures defined in  ${\mathbb   R}^d$. If ${\cal H}
= \{H_i\}_{i=1}^k$ is a collection of $k$ hyperplanes in ${\mathbb
R}^d$ in general position, the connected components of the
complement ${\mathbb  R}^d\setminus\cup~{\cal H}$ are called
(open) $k$-orthants. The definition can be meaningfully extended
to the case of degenerated collections ${\cal H}$, when some of
the $k$-orthants are allowed to be empty.

 A collection ${\cal H}$ is an {\em equipartition}, or more precisely
a $k$-equipartition for ${\cal M}$ if $$\mu_i(O) =
\mu_i(\overline{O\mathstrut}) = {\frac{1}{2^k}}\mu_i({\mathbb
R}^d)$$ for each of the measures $\mu_i\in {\cal M}$ and for each
$k$-orthant $O$ associated to ${\cal H}$.
\end{defin}

\begin{defin}
\label{admissible} A triple $(d,j,k)$ of integers is referred to
as {\em admissible} if for any collection ${\cal M} =
\{\mu_i\}_{i=1}^j$ of $j$ continuous measures in ${\mathbb  R}^d$,
there exists a collection of $k$ hyperplanes ${\cal H} =
\{H_i\}_{i=1}^k$  forming an equipartition for all measures in
${\cal M}$.
\end{defin}

\begin{problem}
\label{problem} The general problem is to characterize the set
${\cal A}$ of all admissible triples. If the emphasis is put on
the ambient Euclidean space ${\mathbb  R}^d$, the equivalent
problem is to determine the smallest dimension $d := \Delta(j,k)$
such that the triple $(d,j,k)$ is admissible.
\end{problem}

\bigskip

One of the main results of the paper is the following theorem
which gives a sufficient condition for a triple $(6m+2, 4m+1,2)$
to be admissible.

\begin{theo}
\label{thm:glavna} Suppose that $(d,j,k) = (6m+2,4m+1,2)$ where
$m$ is a positive integer. Then there exists an equipartition of
$j=4m+1$ mass distributions in ${\mathbb  R}^d = {\mathbb
R}^{6m+2}$ if
\[
\Omega(m):= \alpha(2m+1) + 2\beta(2m+1) - \gamma(2m+1)
\]
is not divisible by $4$, where $\alpha, \beta, \gamma$ are the
combinatorial functions counting special classes of cyclic, signed
$\{A, B\}$-words, introduced in Definition~\ref{def:funkc}.
\end{theo}

It turns out (Theorem~\ref{thm:852}) that $\Omega(1) \equiv 2$
(modulo $4$) which implies that $\Delta(5,2) = 8$ and answers a
well known open question, \cite{Ram}.

These results are obtained by topological methods involving
careful analysis of normal data of the associated singular
(weighted) $1$-manifolds, see Section~\ref{sec:algorithm} and
Figures~\ref{fig:matrica1} and \ref{fig:matrica2}. However we
emphasize that the computation of relevant cohomology/bordism
obstruction classes is closely related to a question of
enumerating classes of immersed curves in $\mathbb{R}^2$ with a
prescribed type and number of intersections with the coordinate
axes, Figures~\ref{fig:metamorf} to \ref{fig:aaaa1}, and in turn
to a problem of enumerating classes of signed, cyclic $A,
B$-words, see Example~\ref{exam:circles}, equations
(\ref{eq:8521})--(\ref{eq:8523}) etc.

Along these lines one obtains other, similar equipartition
results, see Proposition~\ref{prop:case0},
Proposition~\ref{prop:recporec} and Corollary~\ref{cor:general2}.
It is worth mentioning that some of these results have alternative
proofs based on the ideal-valued, cohomological index theory. This
approach  has a merit that it produces fairly general and in some
case quite accurate upper bounds for the function $\Delta(j,k)$,
Theorem~\ref{thm:index}.

\subsection{History of the Problem}

The general problem of studying equipartitions of masses by
hyperplanes, or in our reformulation the problem of  determining
the function $d = \Delta(j,k)$, was formulated by Branko Gr\"
unbaum in \cite{Gru}. Hugo Hadwiger proved that $\Delta(2,2) = 3$,
\cite{Hadw}, which also implies $\Delta(1,3)=3$.  The case $k=1$
is answered by the ``ham sandwich theorem'', \cite{Bo-Ul}, which
says that $\Delta(d,1) = d$. Edgar Ramos, \cite{Ram}, building on
the previous special results in \cite{Yao}, introduced new ideas
and considerably advanced our knowledge about the function $d =
\Delta(j,k)$. He showed for example that $\Delta(1,4) \leq 5,\,
\Delta(5,2) \leq 9,\, \Delta(3,3) \leq 9$. The {\em moment curve}
considerations, see \cite{Ram} or Section~\ref{moment}, lead to
the following general lower bound
\begin{equation}
\label{upper bound} \Delta(j,k)\geq j(2^k-1)/k .
\end{equation}
As far as the general upper bounds are concerned, Ramos proved
that for $j=2^m$, one has the inequalities
\begin{equation}
\label{ramos}
\begin{array}{cccc}
 3j/2\leq \Delta(j,2) \leq 3j/2 &&& 7j/3\leq \Delta(j,3) \leq
 5j/2\\
15j/4\leq \Delta(j,4) \leq 9j/2 &&& 31j/5\leq \Delta(j,5) \leq
15j/2
\end{array}
\end{equation}
It was conjectured in \cite{Ram} that the bound given by
(\ref{upper bound}) is tight.

In computational geometry, the equipartition problem arose in
relation to the problem of designing efficient algorithms for
half-space range queries, \cite{Wil}. In the meantime better
partitioning techniques were found, \cite{Yao-Yao}, \cite{Mat1},
\cite{Mat2}, \cite{Mat3} but the complexity of the equipartition
problem remains of great theoretical interest.

In combinatorial geometry, problems of partitions of sets of
points and dissections of mass distributions have a long tradition
and occupy one of central positions in this field. In this
category are R.~Rado's ``theorem for general measures''
\cite{Rad}, B.~Gr\" unbaum's ``center point theorem'' for convex
bodies \cite{Gru}, etc. Survey articles \cite{Eck}, \cite{Mat},
\cite{Ziv1}, \cite{Ziv2}, \cite{Ziv3} as well as the original
papers \cite{Bar-Mat}, \cite{Ram}, \cite{VZ}, \cite{ZV1}, cover
different aspects of the equipartition problem and give a good
picture of some of more recent developments in this area. For
example \cite{Bar-Mat} deals with partitions by $k$-fans with
prescribed measure ratios, \cite{VZ} studies equipartitions by
regular wedge-like cones and the relation with the well known
Knaster's conjecture, \cite{VZ1} deals with general conical
partitions, the ``center transversal theorem'' proved in
\cite{ZV1} reveals a hidden connection between the center point
theorem and the ham sandwich theorem, etc. Most of these results
are obtained by topological methods. This orientation ultimately
led in \cite{Ziv4}, see also \cite{Mak2002}, to the formulation of
a program which unifies much of the discrete and continuous theory
in the context of so called ``combinatorial geometry on vector
bundles''.

\subsection{Mass Distributions}

\bigskip
A continuous mass distribution, referred to in
Problem~\ref{problem} and Theorem~\ref{thm:glavna}, is a finite
Borel measure $\mu$ defined by the formula $\mu(A) = \int_A f \,
d\mu$ for an integrable density function $f : {\mathbb
R}^d\rightarrow R$. More generally, a mass distribution can be a
Borel measure $\nu$ on ${\mathbb  R}^d$ which is a {\em weak
limit}, \cite{Bill}, of a sequence $\nu_n$ of continuous measures
in the sense that $$ \lim_{n\to\infty}\int_{{\mathbb  R}^d} f\,
d\nu_n = \int_{{\mathbb  R}^d} f\, d\mu $$ for every bounded
continuous function $f : {\mathbb  R}^d\rightarrow R$. From here
one easily deduces the inequalities, $\limsup\, \nu_n(F) \leq
\nu(F)$ for each closed set, and $\liminf\, \nu_n(U) \geq \nu(U)$
for each open set in ${\mathbb  R}^d$.

Given a continuous measure/mass distribution $\mu$, a collection
of $k$ hyperplanes is, according to Definition~\ref{equi}, an
equipartition of $\mu$ if each $k$-orthant has the fraction
$1/2^k$ of the total mass of $\mu$. More generally, one can define
an equipartition for any measure $\nu$ if the inequality $\nu(O)
\leq (1/2^d)\mu({\mathbb  R}^d)$ is valid for each of $2^k$ open
orthants determined by the collection of $k$-hyperplanes. This
shows that Theorem~\ref{thm:glavna} and other equipartition
results which are valid for continuous measures can be easily
extended to more general mass distributions. For example the
statement $\Delta(7,2)=11$ implies that for any set $A =
\{a_{ij}\}_{(i,j)\in [4]\times[7]}\subset {\mathbb  R}^{11}$, a
``matrix'' of points in ${\mathbb  R}^{11}$, there exist two
hyperplanes $H_1$ and $H_2$ such that each of the associated open
quadrants contains at most one point from each of the $7$ columns
of the matrix $A$.

\medskip
In this paper we restrict our attention to measures which are weak
limits of continuous measures which appears to be sufficient for
all combinatorial applications. Among the interesting examples are
measurable sets, counting measures, $k$-dimensional Hausdorff
measures etc.

\section{Proof Technique}
\label{sec:technique}

In this section we collect most of the results needed for the
proof of Theorem~\ref{thm:glavna} and other related equipartition
results. For the readers convenience, in the first subsection we
outline the general proof scheme while other subsections can be
seen as an elaboration of some of the ideas used in individual
steps.

\subsection{The General Scheme of the Proof}
\label{scheme}

Here is a general proof scheme that we follow in this paper and
which in principle can be, with necessary modifications, applied
to many problems about (equi)partitions of masses.

\begin{enumerate}

\item
The equipartition problem, Problem~\ref{problem}, is reduced to
the question of (non)existence of an equivariant map or
equivalently, to the question of the (non)existence of a nowhere
zero, continuous cross-section of a vector bundle. This is a
topological problem.

\item
The topological problem is reduced to the question of computing
relevant (co)homology, or $G$-bordism obstruction classes.

\item The obstruction classes belong to the associated (co)homology
or bordism groups. The computation of these groups usually
involves a mixture of topological and algebraic ideas.

\item
The obstruction class can be, under mild conditions,
computed/identified by a careful analysis of the equipartition
problem for a special, sufficiently ``generic'' collection of
measures. Our primary choice  are uniform/interval measures
distributed along the moment curve $\Gamma_d = \{t,t^2,\ldots ,t^n
\mid t\in R\}\subset {\mathbb  R}^d$. A geometric problem arises,
involving the analysis of possible ways a curve, carrying
collections of intervals, can be immersed in the coordinate two
plane.

\item
The analysis of the solution set of all equipartitions for the
special choice of measures, linked with the problem of immersions
of curves in the previous step,  leads to a problem of enumeration
of classes of circular words in alphabet $\{A,B\}$. This is a
problem of enumerative combinatorics typically solved by Poly\' a
enumeration or M\" obius inversion technique.

\item
The relevant obstruction class is linked in Step~(5) with a
function which has a clear arithmetical/combinatorial meaning
which eventually leads to Theorem~\ref{thm:glavna}.

\end{enumerate}

\subsection{An Approach Based on ${\mathbf{G}}$-bordism} \label{g-bordism}

Here is an example how one can approach an equipartition problem,
in a direct and geometrically transparent fashion, via equivariant
bordism. This is essentially the approach of this paper up to some
technical refinements or detours involving equivariant maps and
equivariant cohomology.

\bigskip
Suppose we want to prove that $\Delta(1,2)=2$ i.e. that for each
measurable set $A\subset {\mathbb  R}^2$, there exist two lines
$L_1$ and $L_2$ which form a ``coordinate system'' so that each
quadrant contains a quarter of the measure of $A$. The problem
itself is of course quite elementary and we use it here merely to
demonstrate the general ideas.

The {\em configuration space} of all candidates for the solution
is the space of all ordered pairs of oriented lines in ${\mathbb
R}^2$. If we fix an embedding ${\mathbb  R}^2\hookrightarrow
{\mathbb  R}^3$ such that ${\mathbb  R}^2$ does not contain the
origin in ${\mathbb  R}^3$, an oriented line $L$ in ${\mathbb
R}^2$ is seen as an intersection of ${\mathbb  R}^2$ with a unique
oriented, central plane $P$ in ${\mathbb  R}^3$. So the
compactified configuration space is the manifold $M = S^2\times
S^2$ of all pairs of oriented, central $2$-planes in ${\mathbb
R}^3$. We note that dihedral group $G:= {\mathbb D}_8$ arises
naturally in this context as a group of all symmetries of a pair
of oriented lines (planes).

For a {\em generic} measurable set $A$, the collection of all
pairs $(L_1,L_2)\in M$ which form an equipartition for $A$ is a
$1$-dimensional $G$-manifold $M_A$. For example if $A$ is a unit
disc $D$, the solution set $M_D$ is a union of $4$ circles. Here
we do not make precise what is meant by a generic measure. Instead
we naively assume, for the sake of this example, that there exists
such a notion of genericity for measurable sets/measures so that
each measurable set $A$ can be well approximated by generic
measures. Moreover, we assume that for any two measurable sets $A$
and $B$ there exists a path of generic measures $\mu_t, \, t\in
[0,1]$, so that $\mu_0$ is an approximation of $A$, $\mu_1$ is an
approximation for $B$ and the solution set
\begin{equation}
\label{} M_{\{\mu_t\}_{t\in[0,1]}}:=\{(L_1,L_2;t)\mid (L_1,L_2)
\mbox{ {\rm is an equipartition for } } \mu_t \} \subset S^2\times
S^2 \times [0,1]
\end{equation}
is a $2$-dimensional manifold (bordism) connecting solution sets
for measures $\mu_0$ and $\mu_1$. The group $\Omega_1({\mathbb
D}_8)$ of classes of $1$-dimensional, free ${\mathbb
D}_8$-manifolds is found to be isomorphic to $\mathbb{Z}/2\oplus
\mathbb{Z}/2$, see Section~\ref{geometric}, and the ${\mathbb
D}_8$-solution manifold $M_D$, associated to the unit disc in
${\mathbb  R}^2$, is shown to represent a nontrivial element in
this group.

It immediately follows than for any measurable set $A\subset
{\mathbb  R}^2$, the solution set $M_A$ is nonempty. Indeed,
suppose $M_A =\emptyset$. Let $\{\mu_t\}_{t\in[0,1]}$ be a path of
generic measures such that $\mu_0$ approximates $A$ and $\mu_1$
approximates $D$. If these approximations are sufficiently good,
we deduce that the solution set $M_{\mu_0}$ is empty and that
$[M_{\mu_1}]$ and $[M_D]$ represent the same element in
$\Omega_1({\mathbb D}_8)$. This is a contradiction since
$M_{\mu_1} =
\partial(M)$ where $M = M_{\{\mu_t\}_{t\in[0,1]}}$, i.e.
$M_D$ would represent a trivial element in $\Omega_1({\mathbb
D}_8)$.

\bigskip
\begin{rem}
\label{rem1} {\rm It is worth noting that the scheme outlined
above, if applicable, shows that a general equipartition problem
can be solved by a careful analysis of the solution set of a well
chosen, particular measure/measurable set, the unit disc $D$ in
our example above. Unfortunately, in higher dimensions the unit
balls do not represent generic measures, i.e. their solution
manifolds are very special and cannot be used for the evaluation
of the relevant obstruction elements in $\Omega_{\ast}({\mathbb
D}_8)$. Instead one uses uniform, interval measures on the {\em
moment curve}, see Section~\ref{moment}. We note also that
although the program outlined above can actually be carried on, we
will use a more standard and equivalent approach via obstruction
classes and equivariant Poincar\' e duality,
Section~\ref{equivariant}. Nevertheless, the idea of a generic
measure is sufficiently interesting and attractive in itself and
we hope to return to it in a later paper.}
\end{rem}

\subsection{Equivariant Maps} \label{equivariant}

It is often more convenient to work with equivariant maps and
their zero sets rather than with measures and their solution
manifolds. Each measure $\mu$ leads naturally to an equivariant
map $A_{\mu} : M\rightarrow V$, where $M$ is a free $G$-manifold
and $V$ a $G$-representation. Working with $G$-equivariant maps
often yields more general results and it is sometimes technically
more convenient. For example one easily constructs $G$-maps which
are transversal to $G$-submanifolds of $V$, thus avoiding the
question of ``generic'' measures from Section~\ref{g-bordism}.

\subsubsection{The central paradigm}
\label{paradigm}

The {\em configuration space/test map} paradigm, \cite{Ziv2}, is
apparently one of the central ideas for applications of
topological methods in geometric and discrete combinatorics.
Review papers, \cite{Alon88}, \cite{Bar93}, \cite{Bjo91},
\cite{Ziv1}, \cite{Ziv2}, \cite{Ziv3} give a detailed picture of
the genesis of these ideas and emphasize their role in the
solutions of well known combinatorial problems like Kneser's
conjecture (L.~Lovasz, \cite{Lov}), ``splitting necklace problem''
(N.~Alon, \cite{Alon87}, \cite{Alon88}), Colored Tverberg problem
(R.~\v Zivaljevi\' c, S.~Vre\' cica, \cite{ZV2}, \cite{VZ94}) etc.

\medskip
The idea can be outlined as follows. One starts with a
configuration space or manifold $M_{\cal P}$ of all candidates for
the solution of a geometric/combinatorial problem ${\cal P}$. For
example, the equipartition problem ${\cal P}$ led in
Section~\ref{g-bordism} to the configuration space  $M_{\cal
P}\cong S^2\times S^2$ of all pairs of oriented planes in
${\mathbb  R}^3$. The following step is a construction of a test
map $f : M_{\cal P}\rightarrow V_{\cal P}$ which measures how far
is a given candidate configuration from being a solution. More
precisely, there is a subspace $Z$ of the test space $V_{\cal P}$
such that a configuration $C\in M_{\cal P}$ is a solution if and
only if $f(C)\in Z$. The inner symmetries of the problem ${\cal
P}$ typically show up at this stage. This means that there is a
group $G$ of symmetries of $X_{\cal P}$ which acts on $V_{\cal
P}$, such that $Z$ is a $G$-invariant subspace of $V_{\cal P}$,
which turns $f : M_{\cal P}\rightarrow V_{\cal P}$ into an
equivariant map. If a configuration $C$ with the desired property
$f(C)\in Z$ does not exist, then there arises an equivariant map
$f : M_{\cal P}\rightarrow V_{\cal P}\setminus Z$. The final step
is to show by topological methods that such a map does not exist.

\subsubsection{Equipartition problem revisited}
\label{sec:revisited}

Our first choice for the configuration space suitable for the
equipartition problem (Problem~\ref{problem}) is the manifold of
all ordered collections $C = (H_1,\ldots ,H_k)$ of oriented
hyperplanes in ${\mathbb  R}^d$. In order to obtain a compact
manifold, we move one dimension up and embed ${\mathbb  R}^d$ in
${\mathbb  R}^{d+1}$, say as the hyperplane determined by the
equation $x_{d+1} = 1$. Then each oriented hyperplane $H$ in
${\mathbb  R}^d \cong \{x\in {\mathbb  R}^{d+1}\mid x_{d+1}=1\}$
is obtained as an intersection $H = {\mathbb  R}^d\cap H'$, where
$H'$ is a uniquely defined, oriented, $(d+1)$-dimensional subspace
of ${\mathbb  R}^{d+1}$. The oriented subspace $H'$ is determined
by the corresponding orthogonal unit vector $u\in S^d\subset
{\mathbb  R}^{d+1}$, so the natural environment for collections $C
= (H_1,\ldots, H_k)$, and our actual choice for the configuration
space is $M_{\cal P}:= (S^d)^k$. The group which acts on the
configuration manifold $M_{\cal P}$ is the reflection group $W_k
:= (\mathbb{Z}/2)^k\rtimes S_k$. The test space $V = V_{\cal P}$
for a single measure is defined as follows. Let $\mu$ be a measure
defined on ${\mathbb  R}^d$ and let $\mu'$ be the measure induced
on ${\mathbb  R}^{d+1}$ by the embedding ${\mathbb  R}^d
\hookrightarrow {\mathbb  R}^{d+1}$. A $k$-tuple $(u_1,\ldots
,u_k)\in (S^d)^k$ of unit vectors determines a $k$-tuple $C =
(H_1',\ldots ,H_k')$ of oriented $(d+1)$-dimensional subspaces of
${\mathbb  R}^{d+1}$. The $k$-tuple $C$ divides ${\mathbb
R}^{d+1}$ into $2^k$-orthants $\mbox{\rm Ort}_{\beta}$ which are
naturally indexed by 0-1 vectors $\beta = \beta_{C}\in
\mathbb{F}_2^k$. Let $b_{\beta} : M_{\cal P}\rightarrow R$ be the
function defined by $b_{\beta}(C) := \mu'(\mbox{\rm Ort}_{\beta})
= \mu(\mbox{\rm Ort}_{\beta}\cap {\mathbb  R}^d)$. Let $B_{\mu} :
(S^d)^k\rightarrow {\mathbb  R}^{2^k}$ be the function defined by
$B_{\mu}(C) = (b_{\beta}(C))_{\beta\in \mathbb{F}_2^k}$. The
$\mu$-test space $V_k = V_{\cal P} \cong {\mathbb  R}^{2^k}$ has a
natural action of the group $W_k := (\mathbb{Z}/2)^k\rtimes S_k$
such that the map $B_{\mu}$ is $W_k$-equivariant. The
$W_k$-representation $V_k$, restricted to the subgroup
$(\mathbb{Z}/2)^k \hookrightarrow W_k$, reduces to the regular
representation $\mbox{\rm Reg}((\mathbb{Z}/2)^k)$ of the group
$(\mathbb{Z}/2)^k$. The ``zero'' subspace $Z_{\cal P}$ is defined
as the trivial, $1$-dimensional $W_k$-representation $V_k^0$
contained in $V_k$. Let $U_k$ be the complementary $W_k$-
representation, $U_k \cong V_k/V_k^0$ and $A_{\mu} : (S^d)^k
\rightarrow U_k$ the induced, $W_k$-equivariant map. By the
construction we have the following proposition which says that
$A_{\mu}$ is a genuine test map for the $\mu$-equipartition
problem.

\begin{prop}
\label{test} A $k$-tuple $C=(H_1,\ldots ,H_k)\in M_{\cal
P}=(S^d)^k$ of oriented hyperplanes is an equipartition of a
measure $\mu$ defined on ${\mathbb  R}^d$ if and only if
$A_{\mu}(C) = 0$. Similarly, given a collection $\{\mu_1,\ldots
,\mu_j\}$ of $j$ measures in ${\mathbb  R}^d$, a $k$-tuple $C$ is
an equipartition for each of the measures $\mu_i$ if and only if
$A(C) = 0$ where $A : (S^d)^k\rightarrow U_k^{\oplus j}$ is the
$W_k$-equivariant map defined by $A(C) := (A_{\mu_i})_{i=1}^j$.
\end{prop}

This proposition motivates the following problem.

\begin{prob}
\label{problem2} Determine or find a nontrivial estimate for the
minimum integer $d := \Theta(j,k)$ such that there does not exist
an $W_k$-equivariant map
\begin{equation}
A : (S^d)^k\rightarrow S(U_k^{\oplus j})
\end{equation}
where $U_k$ is the $W_k$-representation described above, and
$S(U_k^{\oplus j})$ is the $W_k$-invariant (unit) sphere in the
representation $U_k^{\oplus j}$.
\end{prob}

Clearly Proposition~\ref{test}, which says that the inequality
$\Theta(j,k)\geq \Delta(j,k)$ always holds, provides a tool for
proving equipartition results for measures in ${\mathbb  R}^d$.

\subsubsection{Equivariant Poincar\' e duality} \label{duality}

Once a problem is reduced to the question of (non)existence of
equivariant map, one can use some standard topological tools for
its solution. For example one can use the cohomological index
theory for this purpose, \cite{FH2}, \cite{Do}, \cite{Ziv2},
\cite{Ziv3}. This approach is discussed in
Section~\ref{cohomological index}. In this paper our main tool is
elementary equivariant obstruction theory \cite{tDieck87}, refined
by some basic equivariant bordism, and group homology
calculations.

Suppose that $M^n$ is orientable, $n$-dimensional, free
$G$-manifold and that $V$ is a $m$-dimensional, real
representation of $G$. Then the first obstruction for the
existence of an equivariant map $f : M\rightarrow S(V)$, is a
cohomology class
\[
\omega \in H^{m}(M,\pi_{m-1}(S(V)))
\]
in the appropriate equivariant cohomology group, where
$\pi_k(S(V))$ is seen as $G$-module. The action of $G$ on $M$
induces a $G$-module structure on the group
$H_n(M,\mathbb{Z})\cong \mathbb{Z}$ which is denoted by
$\mathcal{Z}$. The associated homomorphism $\theta : G\rightarrow
\{-1,+1\}$ is called the orientation character. Let $A$ be a
(left) $G$-module. The Poincar\' e duality for equivariant
(co)homology is the following isomorphism, \cite{Wall},

\begin{equation}
\label{poincare} H^k_G(M,A)\stackrel{D}\longrightarrow
H_{n-k}^G(M,A\otimes \mathcal{Z})
\end{equation}

Moreover, if $g : M\rightarrow V$ is a smooth map transversal to
$\{0\}\subset V$, then $V := g^{-1}(0)$ is an oriented
$G$-submanifold of $M$ and the dual $D(\omega)$ of the obstruction
class $\omega$ is represented by the fundamental class $[V]$ of
$V$.

 As an illustration, we compute the coefficient $G$-module $N := A\otimes
\mathcal{Z}$ in the case of interest in this paper. According to
Section~\ref{sec:revisited}, Problem~\ref{problem2}, the
equipartitition problem can be reduced to the question of
existence of an equivariant map

\begin{equation}
\label{ponovo}
 A : (S^d)^k\rightarrow S(U_k^{\oplus j}).
\end{equation}
Although the computation in full generality is not much more
difficult, we restrict our attention to the case $k=2$, the only
case which is systematically studied in this paper. The group $G
\cong W_2$, turns out to be a dihedral group ${\mathbb D}_8$ of
order $8$. We work with the presentation of $G \cong {\mathbb
D}_8$ described in Section~\ref{dihedral}. The equivariant map
(\ref{ponovo}) for $k=2$ has the form $A : S^d\times
S^d\rightarrow U_2^{\oplus j}$. The $3$-dimensional vector space
$U_2$ can be decomposed, as a ${\mathbb D}_8$-representation, into
a direct sum $U_2 \cong E_1\oplus E_2$, where $E_2$ is the
standard, $2$-dimensional ${\mathbb D}_8$-representation, and
$E_1$ a $1$-dimensional representation where $\alpha$ and $\beta$
act non trivially while the action of $\gamma$ is trivial. Then
the ${\mathbb D}_8$-structure on the group $A\otimes
\mathcal{Z}\cong \mathbb{Z}$ can be read off the following table.

 {\renewcommand{\arraystretch}{1.5}
\begin{center}

\begin{tabular}{c|c|c|c}
 & $\alpha$ & $\beta$ & $\gamma$ \\ \hline

 $H_{2d}(S^d\times S^d)$ & \quad $(-1)^{d+1}$ & \quad $(-1)^{d+1}$ &
 \, $(-1)^{d^2}$ \\

 $U_2$ & $+1$ & $+1$ & $-1$ \\

 $U_2^{\oplus j}$ & $+1$ & $+1$ & $(-1)^j$ \\ \hline

$A\otimes \mathcal{Z}$ & \quad $(-1)^{d+1}$ &\quad $(-1)^{d+1}$ &
\quad $(-1)^{d+j}$

\end{tabular}

\end{center}
}  In Section~\ref{sec:case1}, we will be particularly interested
in the case $\Delta = 2d-3j =1$, i.e. in the tuples $(d,j) =
(3m+2,2m+1)$, where $m$ is a non negative integer. Then the last
row of the table describing the ${\mathbb D}_8$-module structure
on the coefficient group $N$ has the form

\begin{center}
\begin{tabular}{c|c|c|c}

$N = A\otimes \mathcal{Z}$ & \quad $(-1)^{m+1}$ \quad & \quad
$(-1)^{m+1}$ \quad & \quad $(-1)^{m+1}$ \quad

\end{tabular}

\end{center}

If $m$ is odd, $N$ is a trivial ${\mathbb D}_8$-module which is
denoted simply by $\mathbb{Z}$. If $m$ is even, then $N$ is a
group isomorphic to $\mathbb{Z}$ while all the generators $\alpha,
\beta, \gamma$ act nontrivially. This module structure is in
Section~\ref{comp} denoted by ${\cal Z}$.

\subsection{Moment Curve} \label{moment}

The moment curve $\Gamma_d = \{(t,t^2,\ldots, t^d) \mid t\in R\}$
and the closely related Carath\' eodory curve $C_n =
\{(\cos{t},\sin{t},\cos{2t},\sin{2t}\ldots,\cos{nt},\sin{nt}) \mid
t\in [0,2\pi]\}$, have numerous applications in geometric
combinatorics, \cite{Christos}, \cite{Ziegler95},
\cite{Ziv-advances}. The key property of these curves is that each
hyperplane $H$ intersects $\Gamma_g$ in at most $d$ points;
respectively $2n$ points in the case of Carath\' eodory curve
$C_n$. Let $I_1, I_2, \ldots , I_j$ be a collection of disjoint
intervals on the moment curve $\Gamma_d, \, I_i = [a_i, b_i], \,
a_1 < b_1 < a_2 < b_2 \ldots a_j < b_j $. Let $\mu_i$ be the
uniform probability measure on $I_i$ and $\hat{\mu}_i$ the induced
measure on ${\mathbb  R}^d$, $\hat{\mu}_i(B) := \mu_i(B\cap I_i)$.
A collection $\mathcal{H} = \{H_1, H_2,\ldots , H_k\}$ of $k$
hyperplanes in ${\mathbb  R}^d$ can have at most $kd$ intersection
points with the curve $\Gamma_d$. If $\mathcal{H}$ is an
equipartition of all measures $\hat{\mu}_i$, then the number of
intersection points is at least $j(2^k - 1)$. It follows that if
an equiaprtition exists then $kd\geq j(2^k-1)$ and we obtain the
following lower bound for the function $d = \Delta(j,k)$.

\begin{prop}\mbox{\rm (\cite{Ram})}
\begin{equation}
\label{lower bound} \Delta(j,k)\geq j(2^k-1)/k
\end{equation}
\end{prop}

\subsection{Circular $\{\mathbf{A},\mathbf{B}\}$-words} \label{ab}

Let ${\cal A}_n$ be the set of all words of length $2n$ in the
alphabet $\{A,B\}$ and let ${\cal R}_n$ be the subset of all words
with the same number of occurrences of letters $A$ and $B$. These
words will be occasionally referred to as {\em balanced} words.

\begin{defin}
\label{conjugation} Given a word $w = x_1x_2\ldots x_{2n}$ in
${\cal A}_n$, let $C(w):= x_2x_3\ldots x_{2n}x_1$ be its cyclic
permutation. The conjugation $\ast : {\cal A}_n\rightarrow {\cal
A}_n$ is inductively defined by $A^{\ast} = B$, $B^{\ast} = A$,
$(u\, v)^{\ast} = u^{\ast}v^{\ast}$, i.e. the conjugation is an
involution on ${\cal A}_n$ which replaces each occurrence of a
letter $A$ in $w$, by a letter $B$ and vice versa. The operator
$C$ is a generator of a $\mathbb{Z}/2n$-action on ${\cal A}_n$
while $C$ and $\ast$ together generate a $G := \mathbb{Z}/2n
\times \mathbb{Z}/2$ action on both ${\cal A}_n$ and the set
${\cal R}_n$ of balanced words.
\end{defin}

\begin{defin}
A {\em circular word} is a word in ${\cal A}_n$ up to a cyclic
permutation. More precisely, a circular word is a
$(\mathbb{Z}/2n)$-orbit in the $(\mathbb{Z}/2n)$-set ${\cal A}_n$.
If $w\in {\cal A}_n$, then the associated circular word in ${\cal
A}_n/(\mathbb{Z}/2n)$ is denoted by $[w]$. Let $R(n) := \vert{\cal
R}_n/(\mathbb{Z}/2n)\vert$ be the number of balanced, circular
words of length $2n$.
\end{defin}

\begin{defin}
If $w\in {\cal A}_n$, let \, $\mbox{\rm Per}_{\mathbb{Z}/2n}(w) :=
\mbox{\rm min}\{\, l\mid C^l(w) = w\}$. The number $p = \mbox{\rm
Per}_{\mathbb{Z}/2n}(w)$ is called the period of $w$ and $w$ is
referred to as a word of primitive period $p$.  A word $w\in {\cal
A}_m$ of primitive period $2m$ is called {\em primitive}. If $w\in
{\cal A}_m$ is primitive, then the associate circular word $[w]$
is also called primitive. Let $P(m)$, respectively $Q(m)$, be the
number of primitive, circular, words of length $m$, respectively
the number of primitive, circular, balanced words of length $2m$.
\end{defin}

Poly{\' a}'s enumeration theory, \cite{Gessel95}, \cite{Krishna},
deals with the problem of enumerating the $G$-orbits of classes of
weighted words/functions. An initial example is the following
formula for the number $R(n)$ of, balanced, circular
$\{A,B\}$-words
\begin{equation}
R(n)  =  \frac{1}{2n}\sum_{m\vert n} {{2m}\choose{m}}\phi(n/m)
\end{equation}
where $\phi$ is Euler totient function. An alternative approach to
the problem of counting circular words is via M\" obius inversion
theorem, \cite{Aigner79}, \cite{Gessel95}, \cite{Stan1}. The set
${\cal R}_m$ of all, balanced words of length $2m$ is a disjoint
union of words of primitive period $2k$ for some divisor $k$ of
$m$. Hence,
\begin{equation}
{{2m}\choose{m}} = \sum_{k\vert m}(2k)Q(k)
\end{equation}
and the M\" obius inversion yields the following equation
\begin{equation}
Q(m)  =  \frac{1}{2m}\sum_{k\vert m} {{2k}\choose{k}} \mu(m/k).
\end{equation}

Similarly,
\begin{equation}
\label{eqn:similar} {2^m} = \sum_{k\vert m} k P(k) \qquad \mbox{
yields } \qquad P(m) = \frac{1}{m}\sum_{k\vert m} {2^k} \mu(m/k)
\end{equation}

Keeping in mind the connection between self-conjugated, balanced,
circular $AB$-words with the bordism obstruction classes $o\in
\Omega_1({\mathbb D}_8)$ established in Section~\ref{sec:results},
we now focus our attention on the action of the involution $\ast$
on the set ${\cal R}_n/(\mathbb{Z}/2n)$ of circular words.

\begin{defin}
\label{special} A word $w\in {\cal R}_n$ is {\em special} if
$C^r(w) = w^{\ast}$ for some $r$. A word $w$ is special iff the
associated circular word $[w]\in {\cal R}_n/(\mathbb{Z}/2n)$ is
{\em self-conjugated} in the sense that $\ast([w]) = [w]$. A
primitive, special word in ${\cal R}_m$ is called
$\ast$-primitive. Let $A(m)$ be the number of all $\ast$-primitive
circular words in ${\cal R}_m$. A circular, special word is also
called self-conjugated.
\end{defin}

\begin{lema}
\label{aa^} A word $w\in {\cal R}_n$ is special if and only if it
has the form $(aa^{\ast})(aa^{\ast})\ldots (aa^{\ast})$ for some
$a$. This representation is referred to as a {\em special
representation} of the word $w$. A word is $\ast$-primitive if it
has a unique special representation of the form $w = b b^{\ast}$.
\end{lema}

\noindent It follows from Lemma~\ref{aa^} that the number of
circular, self-conjugated words in ${\cal R}_n/(\mathbb{Z}/2n)$ of
primitive period $2m$, where $m\vert n$, is also $A(m)$.

\medskip

In light of applications in Section~\ref{sec:results}, it would be
interesting to have an explicit, simple formula for the function
$A(m)$. This problem is certainly amenable to more refined methods
from Poly{\' a} enumeration/M\" obius inversion theory, see
\cite{Stan1} or \cite{Krishna}, Section~II.5. However, for our
purposes it is sufficient to determine the residuum of $A(m)$
modulo $2$.

\begin{prop}
\[
A(m) \equiv  P(2m)  \quad \mbox{\rm modulo } \, 2
\]
\end{prop}

\noindent {\bf Proof:} Let ${\cal A}_n(m)$ be the set of all
circular words in ${\cal A}_n/(\mathbb{Z}/2n)$ of primitive period
$2m$ and let ${\cal A}_n^\ast(m)$ be its subset of self-conjugated
circular words. For the proof of the proposition, it is sufficient
to observe that $\vert{\cal A}_n(m)\vert P(2m)$, $\vert{\cal
A}_n^\ast(m)\vert = A(m)$ and that ${\cal A}_n^\ast(m)$ is the
fixed point set for the involution $\ast : {\cal A}_n(m)
\rightarrow {\cal A}_n(m)$. \hfill $\square$

\begin{prop}
\label{prop:racun} $P(2k)\equiv_2 1$ if either $k$ is odd,
square-free integer or $k=2q$ where $q$ is odd, square-free
integer. Otherwise, $P(2k)\equiv_2 0$.
\end{prop}

\medskip\noindent
{\bf Proof:} The proof follows from the explicit formula for the
function $P(m)$ given in equation (\ref{eqn:similar}) and well
known properties of the M\" obius function. \hfill $\square$

\subsection{Homology and Bordism Computations}
\label{homology-bordism-comp}

By equivariant Poincar\' e duality, Section~\ref{duality}, the
dual $D(\omega)$ of the first obstruction cohomology class
$\omega\in H^m_G(M,\pi_{m-1}S(V))$ lies in the group $H^G_{n-m}(M,
\pi_{m-1}S(V)\otimes\mathcal{Z})$. If $M$ is $(n-m)$-connected,
then there is an isomorphism (\cite{Brown82}, Theorem~{II.5.2})
$$H^G_{n-m}(M, \pi_{m-1}S(V)\otimes\mathcal{Z}) \stackrel{\cong}
\longrightarrow H_{n-m}(G, \pi_{m-1}S(V)\otimes\mathcal{Z}).$$
This allows us to interpret $D(\omega)$ as an element in the
letter group. Moreover, if the coefficient $G$-module
$\pi_{m-1}S(V)\otimes\mathcal{Z}$ is trivial, then the homology
group $H_{n-m}(G,\mathbb{Z})\cong H_{n-m}(BG,\mathbb{Z})$ is for
$n-m\leq 3$ isomorphic, \cite{Conner}, to the oriented $G$-bordism
group $\Omega_{n-m}(G)\cong \Omega_{n-m}(BG)$.

Our objective is to identify the relevant obstruction classes.
Already the algebraically trivial case $H_0(G,M)\cong M_G$, where
$M_G = \mathbb{Z}\otimes_G M$ is the group of coinvariants, may be
combinatorially sufficiently interesting. However, the most
interesting examples explored in this paper involve the
identification of $1$-dimensional obstruction classes. Since these
classes in practice usually arise as the fundamental classes of
zero set manifolds, our first choice will be the bordism group
$\Omega_1(G)$.

\subsubsection{Dihedral group ${\mathbf {\mathbb D}_8}$}
\label{dihedral}

In this section we focus our attention on the dihedral group $G =
W_2 = {\mathbb D}_8$. As usual, the group ${\mathbb D}_8\cong
(\mathbb{Z}/2)^{\oplus 2}\rtimes \mathbb{Z}/2$ is identified to
the group of all symmetries of a square with the generators
$\alpha, \beta$ of $({\mathbb{Z}/2})^{\oplus 2}$ seen as
reflections with respect to the $x$-axes and $y$-axes
respectively, while $\gamma$ is the reflection in the line $x=y$.

Low dimensional homology groups are usually not difficult to
compute, say via Hochschild-Serre spectral sequence. However, in
our applications we need more precise information about the
representing cycles of these groups, and this is the reason why
work directly with chains and resolutions.

A convenient $Z[G]$-resolution of $Z$ for the group $G = {\mathbb
D}_8 = ({\mathbb{Z}/2})^{\oplus 2}\rtimes {\mathbb{Z}/2}$, or more
geometrically a convenient $EG$-space with an economical
$G$-CW-structure, can be described as follows.

Let $EG := S^{\infty}\times S^{\infty} \times S^{\infty}$. The
action of $G$ on $EG$ is described as follows
\begin{equation}
\label{action}
\begin{array}{ccc}
\alpha (x,y,z) & = & (-x,y,z) \\ \beta (x,y,z) & = & (x,-y,z) \\
\gamma (x,y,z) & = & (y,x,-z)
\end{array}
\end{equation}

A presentation of ${\mathbb D}_8$ as the group freely generated by
the generators $\alpha, \beta, \gamma$ subject to the relations
\begin{equation}
\label{presentation}
\begin{array}{ccccccc}
\alpha^2 & = & \beta^2 & = & \gamma^2 & = & 1 \\
         &   & \alpha \beta & = & \beta \alpha &  & \\
         &   & \alpha \gamma & = & \gamma \beta &  & \\
         &   & \beta \gamma & = & \gamma \alpha &  &
\end{array}
\end{equation}
shows that the action described by (\ref{action}) is well defined.
There is a natural $G$-invariant $CW$-structure on $EG =
(S^{\infty})^{\times 3}$ which is described as the product of
usual ${\mathbb{Z}/2}$-invariant $CW$-structures on $S^{\infty}$.
In more details, let

\begin{equation}
\label{basic}
\begin{array}{ccccccccc}
\dots & \stackrel{1-t}\longrightarrow & Z[{\mathbb{Z}/2}]x_2 &
\stackrel{1+t}\longrightarrow & Z[{\mathbb{Z}/2}]x_1 &
\stackrel{1-t}\longrightarrow & Z[{\mathbb{Z}/2}]x_0 &
\longrightarrow & 0
\end{array}
\end{equation}
be the usual ${\mathbb{Z}/2}$-invariant cellular chain complex of
$S^{\infty}$ with one cell $x_i$ in each dimension. Then the
cellular chain complex ${\cal C}_G = (\{C_n\}_{n\geq 0},\partial
)$ for $EG$ can be seen as a tensor product of cellular chain
complexes of individual spheres with $t$ in (\ref{basic}) replaced
in the corresponding sphere by $\alpha, \beta,$ and $\gamma$
respectively, and the corresponding generators/cells $x_i$ are
denoted respectively by $a_i, b_i, c_i$. So, a typical $n$-cell in
${\cal C}_G$ has the form $g (a_i\times b_j\times c_k)$ or in a
more algebraic fashion $g (a_i\otimes b_j\otimes c_k)$, for some
$g\in G$ where $i+j+k = n$. Note that according to (\ref{action}),
the action of $G$ on ${\cal C}_G$, the cellular chain complex of
$EG$, is described by the equalities
\begin{equation}
\begin{array}{ccc}
\label{action2} \alpha (a_i\otimes b_j\otimes c_k) & = &(\alpha
a_i\otimes b_j\otimes c_k)\\ \beta (a_i\otimes b_j\otimes c_k) & =
& (a_i\otimes \beta b_j\otimes c_k)\\ \gamma (a_i\otimes
b_j\otimes c_k) & = & (a_j\otimes b_i\otimes \gamma c_k)

\end{array}
\end{equation}

We are interested in the first homology groups $H_1(G,M)$ where
$M$ is either $Z$ seen as a trivial, right $G$-module, or $M =
{\cal Z}$, where ${\cal Z}\cong Z$ as a group, while the
$G$-action on ${\cal Z}$ is defined by $$ t \alpha  = t \beta  = t
\gamma  = -t$$ where $t\in {\cal Z}$ is a generator.

\subsubsection{Homology computations}
\label{comp}

The group $H_n(G,M)$ is by definition the appropriate homology
group of the following chain complex.

\begin{equation}
\label{chain-complex}
\begin{array}{ccccccccc}
\dots & \stackrel{\partial}\longrightarrow & M\otimes_G C_2 &
\stackrel{\partial}\longrightarrow & M\otimes_G C_1 &
\stackrel{\partial}\longrightarrow & M\otimes_G C_0 &
\stackrel{\partial}\longrightarrow & 0
\end{array}
\end{equation}

In order to simplify the notation, we will often use the
abbreviation $(a_i,b_j,c_k)$ for $t\otimes_G(a_i\otimes b_j\otimes
c_k)\in M\otimes_G C_n$, where $n = i+j+k$.

\bigskip

{\bf Case~1} {\boldmath $M = \mathbb{Z}$}

\bigskip

The following equalities are easily derived from (\ref{action2}).
{\small
\begin{equation}
\begin{array}{l}
\partial(a_1,b_0,c_0)=(a_0 - \alpha a_0,b_0,c_0) =0 \\
\partial(a_0,b_1,c_0)=(a_0,b_0-\beta b_0,c_0) = 0\\
\partial(a_0,b_0,c_1)=(a_0,b_0,c_0-\gamma c_0)=0
\end{array}
\end{equation}

\begin{equation}
\begin{array}{l}
\partial(a_2,b_0,c_0)=2 (a_1,b_0,c_0) \\
\partial(a_0,b_2,c_0)=2 (a_0,b_1,c_0)\\
\partial(a_0,b_0,c_2)=2 (a_0,b_0,c_1)
\end{array}
\end{equation}

\begin{equation}
\begin{array}{l}
\partial(a_1,b_1,c_0) = (\partial a_1,b_1,c_0) - (a_1,\partial
b_1,c_0) =0 \\
\partial(a_1,b_0,c_1)= -(a_1,b_0,c_0) + (a_0,b_1,c_0)\\
\partial(a_0,b_1,c_1)= -(a_1,b_0,c_0) + (a_0,b_1,c_0)
\end{array}
\end{equation}
}

From these calculations we easily reach the following conclusion.

\begin{prop}
\label{first-prop} $$ H_1(({\mathbb{Z}/2})^{\oplus 2}\rtimes
{\mathbb{Z}/2}, \mathbb{Z}) \cong {\mathbb{Z}/2}\oplus
{\mathbb{Z}/2}$$ The generators are $X = (a_1,b_0,c_0)$ and
$Y=(a_0,b_0,c_1)$.
\end{prop}

\bigskip
{\bf Case 2} {\boldmath $M = {\cal Z}$}

\bigskip
As before, the calculations are based on equations in
(\ref{action2}).

{\small
\begin{equation}
\label{eqn:rel21}
\begin{array}{l}
\partial(a_1,b_0,c_0)=(a_0 - \alpha a_0,b_0,c_0) = 2(a_0,b_0,c_0) \\
\partial(a_0,b_1,c_0)=(a_0,b_0-\beta b_0,c_0) = 2(a_0,b_0,c_0)\\
\partial(a_0,b_0,c_1)=(a_0,b_0,c_0-\gamma c_0)= 2(a_0,b_0,c_0)
\end{array}
\end{equation}

\begin{equation}
\label{eqn:rel22}
\begin{array}{l}
\partial(a_2,b_0,c_0)= (a_1+\alpha a_1,b_0,c_0) = 0 \\
\partial(a_0,b_2,c_0)= (a_0,b_1+\beta b_1,c_0) = 0\\
\partial(a_0,b_0,c_2)=2 (a_0,b_0,c_1 + \gamma c_1) = 0
\end{array}
\end{equation}

\begin{equation}
\label{eqn:rel23}
\begin{array}{ll}
\partial(a_1,b_1,c_0) & = (a_0-\alpha a_0,b_1,c_0) -
(a_1,b_0-\beta b_0,c_0) \\ &= 2(a_0,b_1,c_0) - 2(a_1,b_0,c_0)\\
\partial(a_1,b_0,c_1)& = (a_0-\alpha a_0,b_0,c_1) - (a_1,b_0,c_0-\gamma
c_0)\\ & = 2(a_0,b_0,c_1) -(a_1,b_0,c_0) - (a_0,b_1,c_0)\\
\partial(a_0,b_1,c_1)& = (a_0,b_0-\beta b_0,c_1) - (a_0,b_1,c_0-\gamma
c_0)\\ & = 2(a_0,b_0,c_1) -(a_0,b_1,c_0) - (a_1,b_0,c_0)
\end{array}
\end{equation}
}

As a consequence we have the following proposition.
\medskip

\begin{prop}
\label{second-prop} $$ H_1(({\mathbb{Z}/2})^{\oplus 2}\rtimes
{\mathbb{Z}/2}, {\cal Z}) \cong {\mathbb{Z}/4}$$
\end{prop}

\noindent{\bf Proof:} It follows from (\ref{eqn:rel21}) that the
elements

\begin{equation}
\begin{array}{ccc}
X_{ab} & := & (a_1, b_0, c_0) - (a_0, b_1, c_0) \\
X_{bc} & := & (a_0, b_1, c_0) - (a_0, b_0, c_1) \\
X_{ca} & := & (a_0, b_0, c_1) - (a_1, b_0, c_0)
\end{array}
\end{equation}
generate the group of cycles. Moreover, essentially the only
relation among them is
\begin{equation}
X_{ab} + X_{bc} + X_{ca} = 0 .
\end{equation}
The relations (\ref{eqn:rel23}), expressed in terms of $X_{ab},
X_{bc}, X_{ca}$, read as follows
\begin{equation}
\begin{array}{ccc}
2 X_{ab} = 0 && X_{ca} = X_{bc} .
\end{array}
\end{equation}
It immediately follows that the homology group
$H_1(({\mathbb{Z}/2})^{\oplus 2}\rtimes {\mathbb{Z}/2}, {\cal Z})$
is isomorphic to $\mathbb{Z}/4$ with elements $X_{ca} = X_{bc}$
representing a generator. \hfill $\square$

\subsubsection{Geometric interpretation}
\label{geometric}

\begin{figure}[htb]
\centering
\includegraphics[scale=0.50]{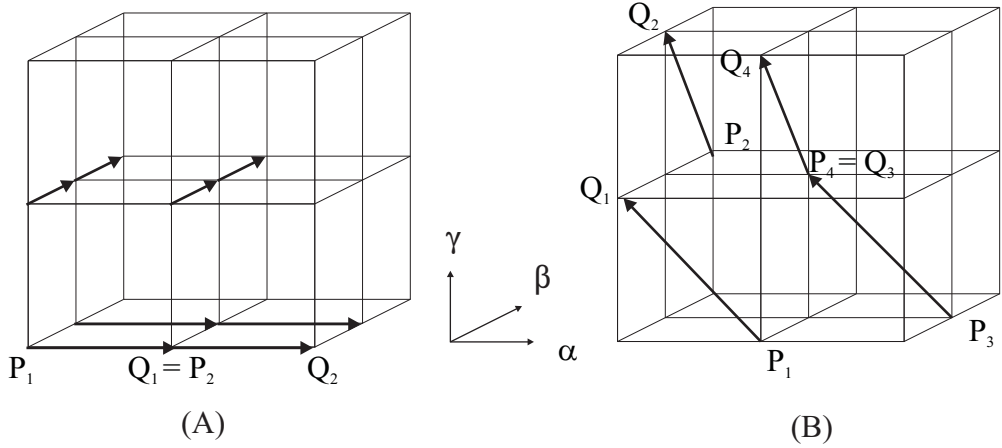}
\caption{Geometric representatives of homology classes}
\label{fig:kocka3}
\end{figure}

The analysis from the previous section allows us to give a precise
geometric interpretation of elements of the groups $H_1({\mathbb
D}_8,\mathbb{Z})\cong \Omega_1({\mathbb D}_8)\cong
\mathbb{Z}/2\oplus \mathbb{Z}/2$ and $H_1(\mathbb{D}_8, {\cal Z})
\cong \mathbb{Z}/4$. By definition $\Omega_1({\mathbb D}_8)$ is
the equivariant bordism group of all $1$-dimensional, oriented,
free ${\mathbb D}_8$-manifolds. Minimal examples are manifolds of
the form ${\mathbb D}_8\times_{\Gamma} S^1 =: M_{\Gamma}$, where
$\Gamma$ is a subgroup of ${\mathbb D}_8$ which acts freely on
$S^1$. The corresponding class in $\Omega_1({\mathbb D}_8)$ is
denoted by $[M_{\Gamma}]$. $\Gamma$ is either the trivial group, a
cyclic group of order $2$ or a cyclic group of order $4$. In order
to associate these manifolds to the elements $X,Y$ and $Z:= X+Y$
of the group $H_1({\mathbb D}_8,\mathbb{Z})$, described in
Proposition~\ref{first-prop}, let us inspect the diagram depicted
on the Figure~\ref{fig:kocka3}~(A). Both diagrams (A) and (B)
represent a $3$-dimensional torus $T^3 = S^1\times S^1\times S^1$,
viewed as a $\mathbb{D}_8$-subcomplex of the complex
$E\mathbb{D}_8 = S^\infty\times S^\infty\times S^\infty$ described
in Section~\ref{dihedral}. Note that the whole $1$-skeleton and a
part of $2$-skeleton of $E\mathbb{D}_8$ are included in this
complex. This implies that all $1$-cycles should be visible in
this picture. For example in Figure~\ref{fig:kocka3}~(A), the
$1$-cell $(a_1,b_0,c_0)$, associated to the generator $X$, is
represented by the vector $P_1Q_1$. Then $P_1Q_1 + \alpha P_1Q_1 =
P_1Q_1 + P_2Q_2$ is a $1$-chain which determines a circle $C_1$
invariant by the group $\mathbb{Z}/2 = \{1,\alpha\}$. The whole
$\mathbb{D}_8$-orbit of $P_1Q_1$  (respectively $C_1$) is in
Figure~\ref{fig:kocka3}~(A) represented by the union of four
circles easily identified as the manifold $M_{(\alpha)}$ where
$(\alpha)$ is the subgroup of ${\mathbb D}_8$ generated by
$\alpha$. A similar analysis shows that $Y = (a_0,b_0,c_1)$ is
represented by $M_{(\gamma)}$. Note that conjugated groups
$\Gamma_1$ and $\Gamma_2$ determine manifolds $M_{\Gamma_1}$ and
$M_{\Gamma_2}$ which are isomorphic as ${\mathbb D}_8$-manifolds
and represent the same elements in $\Omega_1({\mathbb D}_8)$.

The following proposition summarizes the results about the
geometric representatives of elements of the group
$H_1(\mathbb{D}_8, \mathbb{Z})$.

\begin{prop}
\label{recognition1} The nontrivial elements in $H_1(\mathbb{D}_8,
\mathbb{Z})\cong\Omega_1({\mathbb D}_8)\cong
\mathbb{Z}/2\oplus\mathbb{Z}/2$ are
\begin{equation}
X = [M_{(\alpha)}] = [M_{(\beta)}], \quad  Y = [M_{(\gamma)}] =
[M_{(\alpha\beta\gamma)}], \quad Z = [M_{(\alpha\gamma)}] =
[M_{(\beta\gamma)}]
\end{equation}
while the  trivial element in $\Omega_1({\mathbb D}_8)$ is
represented by $1$-manifolds $M_{(e)}$ and $M_{(\alpha\beta)}$.
\end{prop}

\noindent{\bf Proof:} As in the special case of the generator $X$,
one starts with a ``vector description'' $P_1Q_1$ of the
corresponding chain in the torus $T^3$,
Figure~\ref{fig:kocka3}~(A). The smallest $\mathbb{D}_8$-invariant
chain containing $P_1Q_1$ yields the associated manifold
$M_\Gamma$.
 \hfill$\square$

A similar geometric interpretation is possible for the elements of
the group $H_1(\mathbb{D}_8, {\cal Z}) \cong {\mathbb{Z}/4}$. This
time however we have to be more careful. As a motivation we should
keep in mind that the representatives of this group appear ``in
nature'' as {\em weighted} $\mathbb{D}_8$-invariant $1$-manifolds
$N$ where the {weights} take into account the orientation of the
normal bundle $\nu(N)$, seen as a small tubular neighborhood of
$N$ in $M = (S^d)^k$. Recall that we already met weighted chains
$t\otimes_G(a_i\otimes b_j\otimes c_k)$ in Section~\ref{comp}. Let
$M$ be a free, $\mathbb{D}_8$-invariant, oriented, $1$-manifold $M
= M_\Gamma = \mathbb{D}_8\times_\Gamma S^1 = \cup_{g\in G/\Gamma}
S^1_{g}$. In other words $M$ is represented as a disjoint union of
oriented circles $S^1_{g}$ where the indices $g$ run over the
representatives of the corresponding cosets. The circle $S^1 :=
S^1_{e}$ is $\Gamma$-invariant and the circle $S^1_g$ corresponds
to $g\times_\Gamma S^1$ in the representation $M_\Gamma =
\mathbb{D}_8\times_\Gamma S^1$.

\begin{defin}
\label{def:weight} A weighted manifold $M = M_{\Gamma, u}$ is a
formal disjoint sum
 \[ \bigcup_{g\in G/\Gamma}~u_g\otimes S^1_g\]
where the union is taken over the representatives in the cosets
and $u_g\in {\cal Z}$ for each $g\in G$. Moreover, if the function
$u : G\rightarrow {\cal Z}$ is defined so that $u_g = u_h$ if $g$
and $h$ are in the same coset, then there is a compatibility
condition $u_{gh} = g\cdot u_h = gh\cdot u_{e}$ for each $g,h\in
G$. Usually $u_g\in \{t, -t\}$ where $t$ is a preferred generator
of the coefficient module ${\cal Z}$.
\end{defin}

\begin{exam}
\label{exam:weight} {\rm Suppose that $\Gamma = \{1, \gamma\alpha,
\alpha\beta, \gamma\beta\}\cong \mathbb{Z}/4$. Let $u_g := t$ for
each $g\in \Gamma$ and $u_g := -t$ otherwise. Then the weights
$u_g$ obviously satisfy the compatibility condition and
$M_{\Gamma, u} = t\otimes S^1 \cup (-t)\otimes S^1_\alpha$ is a
weighted manifold. }
\end{exam}

\begin{defin}
\label{def:tezine}  A weighted manifold $M = M_{\Gamma, u}$
determines a well defined element $[M]$ in the group
$H_1(\mathbb{D}_8, {\cal Z})$. If $u_g\in\{t, -t\}$ then $M$ is
called a {\em special weighted manifold}. In this case $[M]$ can
be understood as a ``fundamental class'' of $M$ which depends both
on the orientation of $M$ and the weight distribution $u : G
\rightarrow {\cal Z}$. We say that $C\in H_1(\mathbb{D}_8, {\cal
Z})$ is represented by a weighted manifold $M = M_{\Gamma, u}$ if
$C = [M]$.
\end{defin}

\begin{prop}
\label{recognition2} The element $X_{ab}\in H_1(\mathbb{D}_8,
{\cal Z})$ is represented by a special weighted manifold
$M_{(\alpha\beta),u}$. The elements $\pm X_{ca}$ are both
represented by special weighted manifolds of the form
$M_{(\gamma\alpha), u} = u(e)\otimes S^1_e \cup u(\alpha)\otimes
S^1_\alpha$. If we choose a preferred generator $t\in {\cal Z}$
and the orientation on $S^1_e$ so that the action of
$\gamma\alpha$ on $S^1_e$ is a rotation through the angle of
$90^\circ$ in the positive direction, then $X_{ca}$ is represented
by the weighted manifold with the weights  $u_e = t$ and $u_\alpha
= -t$ (Example~\ref{exam:weight}). A representation of $-X_{ca}$
is obtained from the representation of $X_{ca}$ by either changing
the orientation of $S^1_e$ or by using the weights $u_e = -t$ and
$u_\alpha =t$.
\end{prop}

\medskip\noindent
{\bf Proof:} The proof is similar to the proof of
Proposition~\ref{recognition1}. We start with the observation that
one can associate to the generator $X_{ca}$ the vector $P_1Q_1$,
in the usual model of the $\mathbb{D}_8$-space $T^3 = (S^1)^3$,
depicted in Figure~\ref{fig:kocka3}~(B). Then this vector is
completed to a circle $S^1_e := P_1Q_1 + (\gamma\alpha)P_1Q_1 +
(\alpha\beta)P_1Q_1 + (\gamma\beta)P_1Q_1$ which is assumed to be
oriented in the direction of the vector $P_1Q_1$. The proof of the
corresponding statement for the element $X_{ab}$ is similar so we
omit the details. \hfill$\square$

\begin{rem}
\label{rem:vazno} Note that the special weighted manifold
$M_{(\alpha\beta),u}$, representing the element $X_{ab}$, is a
union of four circles. Since $X_{ab}$ is an element of order $2$,
contrary to the case of the generator $X_{ca}$, one does not have
to worry about either the orientation of $S^1_e$ or the precise
value of $u_e\in\{t, -t\}$.
\end{rem}

\subsubsection{The algorithm}
\label{sec:algorithm}

The weighted manifolds arise in applications as follows. Suppose
that $f : N \rightarrow V$ is a $\mathbb{D}_8$-equivariant map,
where $N$ is an oriented, free $\mathbb{D}_8$-manifold and $V$ a
real, linear representation of the group $\mathbb{D}_8$, such that
${\rm dim}(N) - {\rm dim}(V) =1$. Moreover, we assume that the
coefficient module for the obstruction homology classes is the
$\mathbb{D}_8$-module ${\cal Z}$, described in
Section~\ref{duality}.

Assume that $f$ is transverse to $0\in V$, which implies that
$C:=f^{-1}(0)$ is a $\mathbb{D}_8$-invariant, $1$-dimensional
submanifold of $N$. Suppose that $T$ and $t$ are orientations on
$N$ and $V$ respectively. Then there is an induced orientation
$o_C$ on each of the circles $C$, connected components of $M$.
Note that the way how these orientations are transformed under the
action of $\mathbb{D}_8$ is governed by the module ${\cal Z}$, for
example if $g\in\{\alpha, \beta, \gamma\}$ then $g(o_C) = -
o_{g(C)}$. The manifold $M$ can be decomposed into the minimal,
$\mathbb{D}_8$-invariant components, $M = M_1\cup\ldots \cup M_n$.
Then each of the manifolds $M_i$ has the corresponding fundamental
class $[M_i]$, which can be evaluated in the group
$H_1(\mathbb{D}_8, {\cal Z})$ by the following rule. If $M_i$ has
$8$ connected components then $[M_i]$ is equal to $0$. If $M_i$
has $4$ connected components then $[M_i]$ is equal to $X_{ab}$. If
$M_i$ has $2$ connected components, let $C$ be one of them,
compare the orientation $o_C$ with the orientation determined by
the requirement that $\gamma\alpha$ is a rotation through
$90^\circ$ in the positive direction. If these two orientations
agree, then $[M_i] = X_{ca}$, otherwise $[M_i] = - X_{ca}$.

\medskip
In practise we are interested mainly in the fact whether the
obstruction element $X = [M_1]+\ldots +[M_n]$ is different from
zero. Hence, minor changes in the algorithm above, for example
choosing the rotation through $270^\circ$ instead of $90^\circ$ as
we did in Section~\ref{sec:852}, does not affect the answer.

\section{Results And Proofs}
\label{sec:results}

In this section we give a detailed analysis of the case $(d,j,2)$,
i.e. we analyze when there exist in a $d$-dimensional space,
simultaneous equipartitions of $j$-measures by two hyperplanes
$H_1, H_2$. In light of the necessary condition $2d\geq 3j$, it
would be desirable to test pairs $(d,j)$ for which the difference
$\Delta := 2d-3j$ is as small as possible.

\subsection{The case $\Delta = 0$}
\label{sec:case0}

We want to determine for which integers $m$, the triple
$(3m,2m,2)$ is admissible. Following the general scheme outlined
in Section~\ref{scheme}, we choose $2m$ disjoint intervals
$I_1,I_2,\ldots ,I_{2m}$ on the moment curve $\Gamma =
\{(t,t^2,\ldots , t^{3m}\mid t\in R\}\subset \mathbb{R}^{3m}$,
representing the measures $\mu_1,\ldots ,\mu_{2m}$, where
$\mu_i(A)$ is by definition the $1$-dimensional (Lebesgue) measure
of the measurable set $A\cap I_i$.

\begin{figure}[htb]
\centering
\includegraphics[scale=0.60]{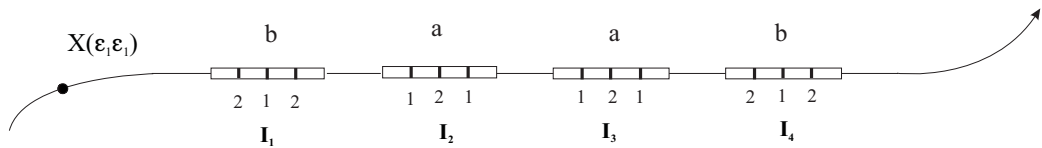} \caption{An equipartition of type
$(\epsilon_1,\epsilon_2)baab$}
\label{fig:linija1}
\end{figure}

Since a hyperplane in ${\mathbb  R}^{3m}$ can intersect the curve
$\Gamma$ in at most $(3m)$-points, we observe that all
intersection points have to be used for partitioning the intervals
$I_j, \, j=1,\ldots, 2m$, an example for $m=2$ is shown in
Figure~\ref{fig:linija1}. Let $(H_1,H_2)$ be a pair of oriented
hyperplanes which is a solution to the equipartition problem for
these intervals.  Then the {\em type} $\tau(H_1,H_2)$ of the
solution $(H_1,H_2)$ is a word in alphabet $\{a,b\}$ and a sign
vector $(\epsilon_1,\epsilon_2), \, \epsilon_i\in \{+,-\}$, which
together encode all the data about the equipartition. For example
in Figure~\ref{fig:linija1}, the first hyperplane $H_1$ intersects
$\Gamma$ in $6$ points which are all labelled by $1$, while the
points in common to $\Gamma$ and $H_2$ are labelled by $2$.

\bigskip

\bigskip

Suppose that $\alpha_i, \, i=1,2$ is the unit vector orthogonal to
$H_i$ pointing in the positive direction. If $X\in \Gamma$ is a
test point preceding all intervals $I_j$, then the type of
$(H_1,H_2)$ is $$ \tau(H_1,H_2) = (\epsilon_1 \epsilon_2) b a a
b$$ where $\epsilon_i := {\rm sign}(\langle \alpha_i,X\rangle)$
and the letters $a$ and $b$ record the way each individual
interval is partitioned. Obviously, in the case $(3m,2m,2)$, any
sign vector and any $ab$-word of length $2m$ with the same number
of occurrences  of letters $a$ and $b$, serves as the type of a
unique solution $(H_1,H_2)$ . Let us denote by ${\cal T}$ the set
of all types. The equivariant Poincar\' e dual to the obstruction
class lies in the group $H_0(({\mathbb{Z}/2})^{\oplus 2}\rtimes
{\mathbb{Z}/2},M)\cong {\mathbb{Z}/2}$ where $M$ is one of the
modules $\mathbb{Z}, {\cal Z}$ described in the
Section~\ref{duality}.

\bigskip
\noindent {\bf Observation:} The dual of the obstruction class is
a nonzero element in $H_0(({\mathbb{Z}/2})^{\oplus 2}\rtimes
{\mathbb{Z}/2},M)\cong {\mathbb{Z}/2}$ if and only if the number
of $G = ({\mathbb{Z}/2})^{\oplus 2}\rtimes {\mathbb{Z}/2}$ orbits
in the $G$-set ${\cal T}$ is {em odd}.

\bigskip

\begin{prop}
\label{prop:case0} The equivariant map $A : (S^{3m})^2 \rightarrow
S((U_2)^{\oplus 2m})$ (Problem~\ref{problem2}) exists if and only
if $m = 2^q - 1$ for some integer $q$. It follows that a triple
$(3\cdot 2^q -3,2^{q+1}-2,2)$ is admissible, i.e.
$\Delta(2^{q+1}-2,2) = 3\cdot 2^q - 3$.
\end{prop}

\noindent {\bf Proof:} The number of $G$-orbits in ${\cal T}$ is
$\frac 12{2m\choose m} = {{2m-1}\choose{m-1}}$. The result follows
from the well known fact that ${n\choose k}_2 = 1$ iff the binary
representation of $k$ is a subword of the binary representation of
$n$. \hfill $\square$

\subsection{The case $\Delta = 1$}
\label{sec:case1}

 We have already observed in Section~\ref{sec:case0} that our method in
the case $\Delta = 2d-3j=0$ does not provide too many interesting
admissible triples. It turns out that the case $\Delta =1$ is much
more interesting from this point of view, so in this section we
focus our attention on the triples of the form
$(d,j,2)=(3m+2,2m+1,2)$. Following essentially the notation from
Section~\ref{g-bordism}, we denote by ${\cal M} = {\cal
M}[(\mu_{\nu})_{\nu=1}^j]$ the space of all $k$-tuples
$(H_1,\ldots ,H_k)$ of oriented hyperplanes in ${\mathbb  R}^d$
which form an equipartition of the $d$-space with respect to each
of the measures $\mu_\nu$. To be precise, we always work with the
compactified configuration space $M_{\cal P} = (S^d)^k$, so one of
the hyperplanes $H_i$ is allowed to be the hyperplane ``at
infinity'', cf. Section~\ref{sec:revisited}.

In our case $(d,j)=(3m+2,2m+1)$ and we assume that the set
$\mu_1,\ldots ,\mu_j$ of test measures is selected as in the
Sections~\ref{moment} and \ref{sec:case1} as the measures
concentrated and uniformly distributed on disjoint subintervals
$I_0,I_1,\ldots ,I_{2m}$ of the moment curve $\Gamma =
\Gamma_{3m+2}$. Our goal is to describe the set ${\cal M}$ of
solutions to the equipartition problem in the case
$(d,j)=(3m+2,2m+1)$. Moreover, in order to be able to compute the
relevant obstruction class $[{\cal M}]$ in one of the groups
$H_1(\mathbb{D}_8; \mathbb{Z})\cong \Omega_1({\mathbb D}_8)$ or
$H_1(\mathbb{D}_8; {\cal Z})$, we will need a more precise
``inventory'' of classes of these solutions.

\subsubsection{Types and inventory of solutions}

\begin{figure}
\centering
\includegraphics[scale=0.57]{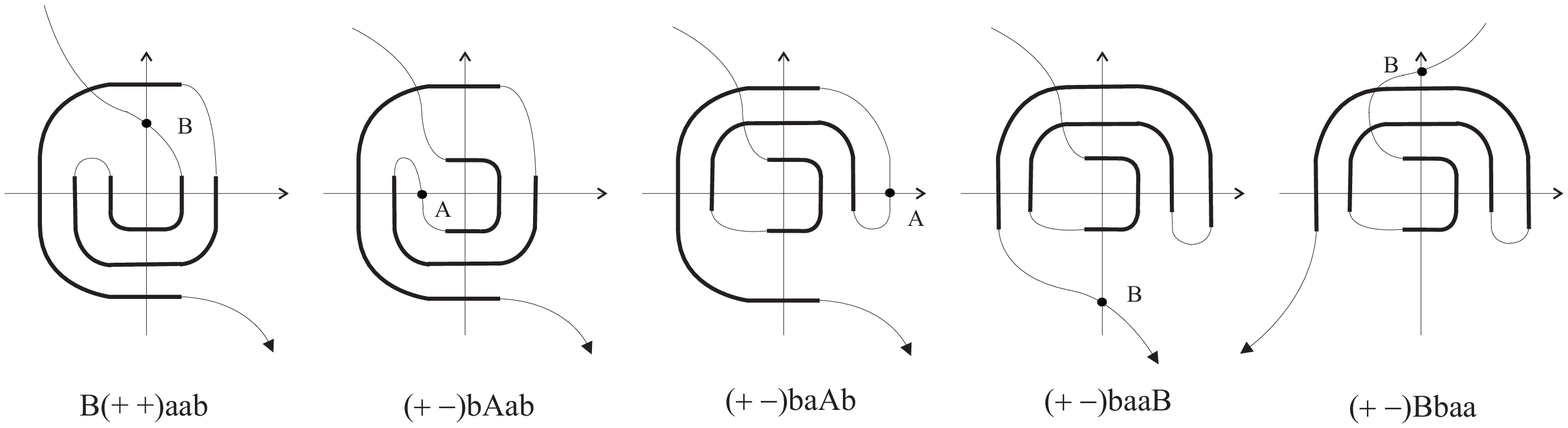} \caption{Metamorphoses of curves}
\label{fig:metamorf}
\end{figure}

 As in Section~\ref{sec:case0}, it is convenient to record the {\em type}
$\tau(H_1,H_2)$ of a solution pair $(H_1,H_2)$ at least in some
typical cases. As opposed to the case $\Delta = 0$, this time
there is one more degree of freedom, arising from the fact that
not all intersection points of hyperplanes $H_1$ and $H_2$ with
the curve $\Gamma$ are necessary for the equipartitions of given
intervals. As a consequence ${\cal M}$ turns out to be a
$1$-dimensional manifold. Figure~\ref{fig:metamorf} serves as a
fairly good sample of ``snapshots '' from a ``movie'' which
describes the solution manifold ${\cal M}$ in the case
$(d,j)=(5,3)$. Each individual drawing in
Figure~\ref{fig:metamorf} represents a single solution. More
precisely, a drawing should be interpreted as an ``orthogonal
image'' of the relevant part of moment curve $\Gamma$ containing
the intervals $I_0, I_1, I_2$ in the $2$-plane $P$ orthogonal to
$H_1\cap H_2$. We emphasize that the hyperplanes $H_1$ and $H_2$,
or rather their intersections with $P$, are represented on these
pictures by the horizontal and vertical coordinate lines, while
the associated orthogonal unit vectors are denoted by $\alpha$ and
$\beta$.

The type $\tau(H_1,H_2)$ of the solution pair $(H_1,H_2)$,
depicted in the first drawing, is by definition $$ \tau(H_1,H_2) =
B(++)aab$$ which takes into account that

\begin{itemize}
\item
the three intervals $I_0,I_1,I_2$ are partitioned according to the
word $aab$ (cf. Section~\ref{sec:case0}),
\item
the initial point of the interval $I_0$ belongs to the first
quadrant and has $(++)$ for the associated sign vector,
\item
the {\em free} intersection point, i.e. the point ($B$ on the
picture) not used for partition of intervals, belongs to $H_2$ and
precedes all intervals $I_i$.
\end{itemize}

The same bookkeeping scheme is used for describing the types of
other solutions shown in Figure~\ref{fig:metamorf}. For example
$\tau(H_1,H_2) = (+-)bAab$ is the type of the solution depicted on
second drawing in Figure~\ref{fig:metamorf}.

\begin{figure}[htb]
\centering
\includegraphics[scale=0.75]{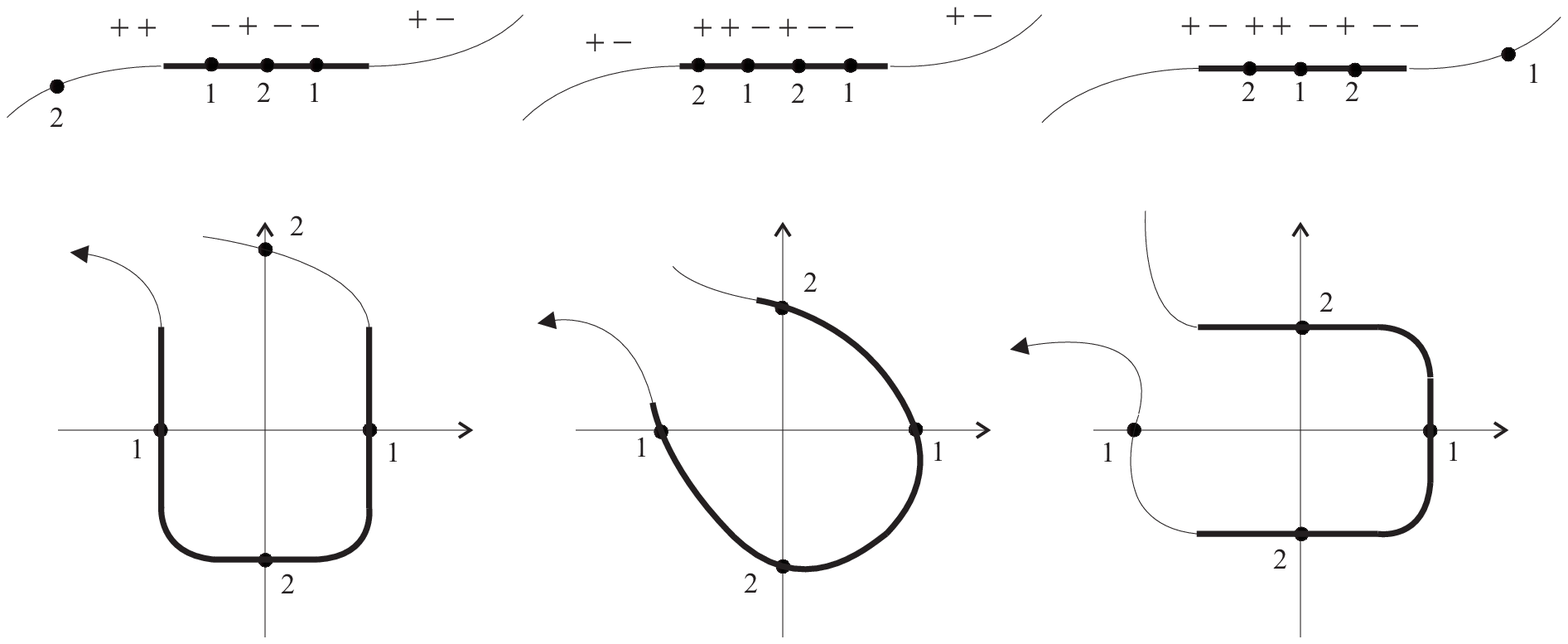} \caption{$B(++)a \leadsto (+-)bA$}
\label{fig:baba1}
\end{figure}

\begin{figure}[bth]
\centering
\includegraphics[scale=0.60]{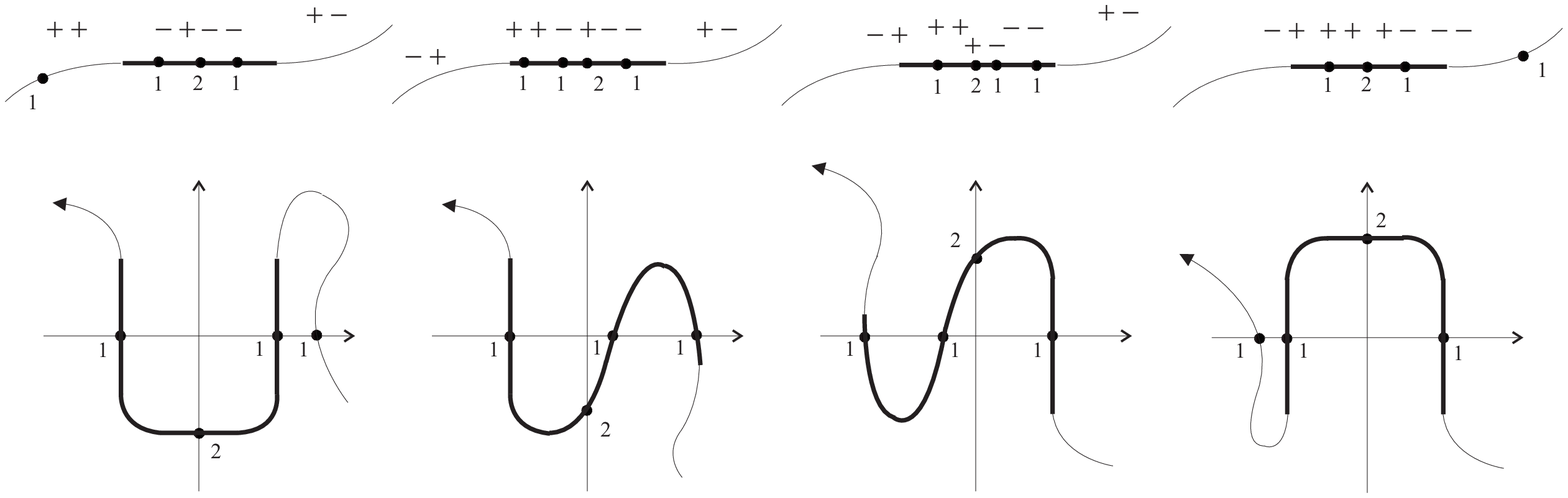} \caption{$A(++)a \leadsto (-+)aA$}
\label{fig:aaaa1}
\end{figure}

Our next observation is that all solutions depicted in
Figure~\ref{fig:metamorf} belong to the same connected component
of the manifold ${\cal M}$. Indeed, these ``metamorphoses'' of
solutions leading from one drawing to another, not visible in
Figure~\ref{fig:metamorf}, are shown in Figures~\ref{fig:baba1}
and \ref{fig:aaaa1}. It turns out that there exist essentially two
types of metamorphoses, the $xy$-type shown in
Figure~\ref{fig:baba1} (marked by the appearance of both letters
$a$ and $b$) and the opposite $xx$-type shown in
Figure~\ref{fig:aaaa1}. The key observation is that the types of
solutions obey a very simple low of transformation. As before, a
{\em free} intersection point is a point $A\in \Gamma\cap H_1$ or
$B\in \Gamma\cap H_2$ which does not belong to any of the
intervals $I_j$.

\begin{prop}
\label{prop-rules} The passage of a {\em free} point through an
interval is recorded on Figures~\ref{fig:baba1}  and
\ref{fig:aaaa1}. The types of associated solutions change
according to the following rules

\begin{equation}
\label{rules}
\begin{array}{ccc}
 B(\epsilon_1,\epsilon_2)a   & \longrightarrow &
 (\epsilon_1,-\epsilon_2)bA \\
 A(\epsilon_1,\epsilon_2)b   & \longrightarrow &
 (-\epsilon_1,\epsilon_2)aB \\
  A(\epsilon_1,\epsilon_2)a   & \longrightarrow &
 (-\epsilon_1,\epsilon_2)aA \\
  B(\epsilon_1,\epsilon_2)b   & \longrightarrow &
 (\epsilon_1,-\epsilon_2)bB
\end{array}
\end{equation}
where $A,B$ are free points, letters $a,b$ record the type of a
partition and $(\epsilon_1,\epsilon_2)\in \{+,-\}^2$ is a sign
vector.
\end{prop}

\noindent {\bf Proof:} The proof is essentially by inspection of
Figures~\ref{fig:baba1} and \ref{fig:aaaa1}. Note that the
transformation rules are equivariant with respect to the group
${\mathbb D}_8$ which allows us to assume that
$(\epsilon_1,\epsilon_2) = (+,+)$. \hfill $\square$

\medskip
\begin{rem}
\label{rem-rules} {\rm Proposition~\ref{prop-rules} records the
change of the type of a solution if the free point passes over the
first interval $I_0$. This is the reason why we see the change of
the sign vector $(\epsilon_1,\epsilon_2)$. Note however, that the
same analysis yields similar rules for the change of types of
solutions when a free point interchanges its position with some of
the intervals $I_i$ for $i\neq 1$
\begin{equation}
\label{rules1}
\begin{array}{cccc}
 Ba    \rightarrow bA &  Ab  \rightarrow aB & Aa \rightarrow aA & Bb
 \rightarrow bB
\end{array}
\end{equation}
}
\end{rem}
These changes do not involve the change of the associated sign
vector. The reader should note that both (\ref{rules}) and
(\ref{rules1}) have the following simple interpretation. They both
say that the underlying word does not change except for a shift of
the capital letter one step to the right.

\medskip

The free point, denoted by a capital letter $A$ or $B$ in the type
$\tau(H_1,H_2)$ of a solution $(H_1,H_2)$, moves to $+\infty$ on
the moment curve $\Gamma$, passes through the infinite point and
appears again on the other side and approaches intervals
$I_1,\ldots , I_j$ from $-\infty$. This passage through an
infinite point is formally justified by the fact that the point
$(0,\ldots ,0,1)\in {\mathbb  R}^{d+1}$ determines the point
$+\infty = -\infty$ which can be seen as an infinite point on the
moment curve $\Gamma$  compactified configuration space $M_{\cal
P} = (S^d)^k$

\begin{exam}
\label{exam-types} {\rm Here is a complete sequence of types of
solutions in the connected component of the solution manifold
${\cal M} = {\cal M}_{(5,3)}$ which contains all solutions
depicted in the Figure~\ref{fig:metamorf}.

\begin{equation}
\label{sequence-of-types}
\begin{array}{ccccccc}
\tau_1=B(++)aab && \tau_5=B(+-)baa && \tau_9=A(++)bba &&
\tau_{13}=A(-+)abb \\

\tau_2=(+-)bAab && \tau_6=(++)bBaa &&\tau_{10}=(-+)aBba &&
\tau_{14}=(++)aAbb \\

\tau_3=(+-)baAb && \tau_7=(++)bbAa && \tau_{11}=(-+)abBa &&
\tau_{15}=(++)aaBb
\\

\tau_4=(+-)baaB && \tau_8=(++)bbaA && \tau_{12}=(-+)abbA &&
\tau_{16}=(++)aabB

\end{array}
\end{equation}

Note that $\tau_{17}=\tau_1 = B(++)aab$ while the types $\tau_1,
\tau_2,\ldots, \tau_5$ correspond to the solutions depicted in
Figure~\ref{fig:metamorf}. The list~(\ref{sequence-of-types}) can
be conveniently abbreviated as follows

\begin{equation}
\label{seq-of-words}
\begin{array}{cccccccccc}
BAAB(+-) &&& BBAA(++) &&& ABBA(-+) &&& AABB(++)
\end{array}
\end{equation}
}
\end{exam}

\begin{defin}
\label{sign-words} A signed word is a word of the form
$w(\epsilon_1,\epsilon_2)$, where $w = x_1x_2\ldots x_{2r}$ is a
balanced $\{A,B\}$-words in the sense of Section~\ref{ab}. Two
signed words, $w_1(\epsilon_1,\epsilon_2)$ and
$w_2(\eta_1,\eta_2)$ are cyclically equivalent if one can be
obtained from the other by subsequent applications of the
following rules

\begin{equation}
\label{rules2}
\begin{array}{cc}
 Aw(\epsilon_1,\epsilon_2) \longleftrightarrow wA(-\epsilon_1,\epsilon_2)
 &
Bw(\epsilon_1,\epsilon_2) \longleftrightarrow
wB(\epsilon_1,-\epsilon_2)
\end{array}
\end{equation}
Cyclic equivalence is an equivalence relation, and an equivalence
class is called a {\em cyclic signed word}. The cyclic signed
word, associated to a signed word $C = w(\epsilon_1,\epsilon_2)$
is denoted by $[C] = [w(\epsilon_1,\epsilon_2)]$.
\end{defin}

Note that the cyclic signed word represents a combinatorial
counterpart of a connected component in a solution manifold. More
precisely we have the following observation.

\begin{observ}
\label{observ1} {\rm  Example~\ref{exam-types} clearly shows that
each circle, representing a connected component of the solution
manifold ${\cal M} = {\cal M}_{d,j}$, can be encoded by an
appropriate {\em cyclic sign word}. Conversely, each cyclic signed
word is associated to some circle in the corresponding solution
manifold. The one to one correspondence between circles in
solution manifolds and cyclic signed words, is ${\mathbb
D}_8$-equivariant in the following sense. The group ${\mathbb
D}_8$, which acts on the solution manifold ${\cal M}$, acts also
both on signed words and on their cyclic equivalence classes. The
action on signed words is described as follows
\begin{equation}
\alpha(w(\epsilon_1,\epsilon_2)) = w(-\epsilon_1,\epsilon_2)
\qquad \beta(w(\epsilon_1,\epsilon_2)) = w(\epsilon_1,-\epsilon_2)
\qquad \gamma(w(\epsilon_1,\epsilon_2)) =
w^{\ast}(\epsilon_2,\epsilon_1)
\end{equation}
where the conjugation $w\mapsto w^{\ast}$ is the operation
described in Definition~\ref{conjugation}. }
\end{observ}

\begin{defin}
\label{def:orbit} The ${\mathbb D}_8$-orbit of a cyclic signed
word is called a {\em generating class of words} or simply a
generating class. The one to one correspondence between cyclic
signed words and circles in solution manifolds, extends to the
correspondence between generating classes and minimal ${\mathbb
D}_8$-invariant submanifolds of the solution manifold. These
minimal ${\mathbb D}_8$-invariant submanifolds have, according to
Propositions~\ref{recognition1} and \ref{recognition2}, the form
$M_{\Gamma}$, for an appropriate subgroup $\Gamma\subset {\mathbb
D}_8$, and they are potential carriers of nontrivial obstruction
classes $X, Y, Z$, respectively $X_{ab}, \pm X_{ca}$.
\end{defin}

\begin{exam}
\label{exam:circles} {\rm The ${\mathbb D}_8$-orbit of the cyclic
signed word {\rm (\ref{seq-of-words})}, i.e. the corresponding
generating class of words, is exhibited in the following table.

\begin{equation}
\label{generating-class}
\begin{array}{ccccccc}
AABB++ && AABB+- && AABB -+ && AABB--\\

ABBA-+ && ABBA-- && ABBA++ && ABBA+-\\

BBAA++ && BBAA+- && BBAA-+ && BBAA-- \\

BAAB+- && BAAB++ && BAAB-- && BAAB -+
\end{array}
\end{equation}
The cyclic signed word (\ref{seq-of-words}) is exhibited as the
first column in the table (\ref{generating-class}). It is easy to
check that its stabilizer is the group $(\gamma)\subset {\mathbb
D}_8$. We infer from here that the corresponding submanifold of
the solution manifold ${\cal M}_{(5,3)}$ is of the type
$M_{(\gamma)}$, in the notation of Proposition~\ref{recognition1}.
Hence, the generating class (\ref{generating-class}) contributes
the element $Y$ to the associated homology or bordism obstruction
class $D(\omega)$ described in
Section~\ref{homology-bordism-comp}. It is easily checked that in
the case $(d,j)=(5,3)$ there is exactly one more generating class
of words shown on the following table,

\begin{equation}
\begin{array}{ccc}
ABAB++ && ABAB+- \\ BABA-+ && BABA-- \\ ABAB-- && ABAB-+ \\ BABA+-
&& BABA++
\end{array}
\end{equation}
The stabilizer of the first column is $(\alpha\beta)$ which,
according to Proposition~\ref{recognition1}, represents a trivial
element in the group $\Omega_1({\mathbb D}_8)$.

The final conclusion is that the obstruction $[{\cal M}] = [{\cal
M}_{(5,3)}] = Y$ is nontrivial which implies that the triple
$(5,3,2)$ is admissible. }
\end{exam}

\subsubsection{The  $(8,5,2)$ case} \label{sec:852}

\begin{theo}
\label{thm:852} The triple $(8,5,2)$ is admissible. In other words
for each collection of $5$ continuous mass distributions $\mu_1,
\ldots , \mu_5$ in $\mathbb{R}^8$ there exist two hyperplanes
$H_1$ and $H_2$ forming an equipartition for each of the measures
$\mu_i$.
\end{theo}

\medskip\noindent
{\bf Proof:} The first obstruction $X$ for the existence of a
$\mathbb{D}_8$-equivariant map $f : S^8\times S^8 \rightarrow
U_2^{\oplus 5}$ is an element of the group $H_1(\mathbb{D}_8,
{\cal Z})$. It turns out that in the $(8,5,2)$-case there exist
precisely three generating classes (Definition~\ref{def:orbit}) of
signed words, displayed on the diagrams (\ref{eq:8521}),
(\ref{eq:8522}) and (\ref{eq:8523}) respectively.

{\footnotesize 
\begin{equation}
\label{eq:8521}
\begin{array}{ccc}
AAABBB++ && AAABBB+- \\  AABBBA-+ && AABBBA-- \\ ABBBAA++ && ABBBAA+- \\
BBBAAA-+ && BBBAAA--\\ BBAAAB-- && BBAAAB-+ \\ BAAABB-+ && BAAABB-- \\
AAABBB-- && AAABBB-+ \\ AABBBA+- && AABBBA++ \\ ABBBAA-- && ABBBAA-+ \\
BBBAAA+- && BBBAAA++ \\ BBAAAB++ && BBAAAB+- \\ BAAABB+- &&
BAAABB++
\end{array}
\end{equation}

\medskip
\begin{equation}
\label{eq:8522}
\begin{array}{ccccccc}
AABABB++ && AABABB+-  &&  AABBAB++ && AABBAB+- \\

 \cdot && \cdot  &&   \cdot && \cdot \\

BABBAA++ && BABBAA+- && BBABAA++ && BBABAA+- \\

 \cdot && \cdot &&  \cdot && \cdot \\

BBAABA-- && BBAABA-+ && ABAABB++ && ABAABB+-\\

 \cdot && \cdot &&  \cdot && \cdot \\

AABABB-- && AABABB-+ && AABBAB-- && AABBAB-+ \\

 \cdot && \cdot &&  \cdot && \cdot \\

BABBAA-- && BABBAA-+ && BBABAA-- && BBABAA-+\\

 \cdot && \cdot &&  \cdot && \cdot \\

BBAABA++ && BBAABA+- && ABAABB-- && ABAABB-+\\

 \cdot && \cdot &&  \cdot && \cdot
\end{array}
\end{equation}

\begin{equation}
\label{eq:8523}
\begin{array}{ccc}
ABABAB++ && ABABAB+- \\ BABABA-+ && BABABA-- \\ ABABAB-- &&
ABABAB-+ \\ BABABA+- && BABABA++
\end{array}
\end{equation}
}

\medskip\noindent
If $M_1, M_2$ and $M_3$ are the associated minimal
$\mathbb{D}_8$-manifolds and $[M_i]$ are the associated weighted
fundamental classes, then $X = [M_1] + [M_2] + [M_3]$. Since $M_2$
has $4$ connected components we conclude that $[M_2] = X_{ab} =
2X_{ca}$. We will show now that $[M_1]$ and $[M_3]$ cancel out
following the procedure described by the algorithm in
Section~\ref{sec:algorithm}.

Let $C_1$ be the circle corresponding to the first column in the
diagram (\ref{eq:8521}), oriented following the lexicographical
order of signed words. Note that the element $\gamma\alpha\in
\mathbb{D}_8$ rotates this circle through the angle of
$270^\circ$. Similarly, $C_2$ is the circle corresponding to the
first column in the diagram (\ref{eq:8523}) oriented by the same
rule. Note that this circle is also rotated by $\gamma\alpha$
through the angle of $270^\circ$.

Following the procedure described in Section~\ref{sec:algorithm},
we should choose orientations $T$ and $t$ on $N = S^8\times S^8$
and on $U_2^{\oplus 5}$. Let us suppose that aside from the
intervals $I_i = [a_i, b_i]$ on the moment curve $\Gamma_8$,
representing the measures $\mu_i$, we choose one more (open)
interval $J_0 = (a_0, b_0)$ such that $b_0 < a_1$. Here for
simplicity we identify $\Gamma_8$ with the associated parameter
space $\mathbb{R}$. Let $J_1 := (a_1, b_5)$ and choose a point $z$
in the interval $(b_0,a_1)$. Define ${\cal A}_1\subset S^8$ as the
space of all oriented hyperplanes $H_1$ such that $z$ is in the
positive halfspace $H_1^+$ and $H_1\cap \Gamma_8$ consists of
precisely $8$ distinct points $x_1< x_2 <\ldots < x_8$ such that
$x_1\in J_0$ and $x_i\in J_1$ for each $i\geq 2$. Similarly, let
${\cal A}_2$ be a collection of all hyperplanes $H_2$ such that
again $z$ is in the positive halfspace $H_2^+$ and the
intersection $H_2\cap \Gamma_8$ consists of $8$ distinct points
$y_1<y_2<\ldots <y_8$ but this time $y_i\in J_1$ for each $i\geq
1$. Both ${\cal A}_1$ and ${\cal A}_2$ are open (convex) cells in
$S^8$ and ${\cal A}_1\times {\cal A}_2$ is a coordinate chart on
$S^8\times S^8$ with coordinates $x_1, \ldots x_8, y_1,\ldots ,
y_8$. The tangent vectors $\partial/\partial x_1,\ldots
,\partial/\partial x_8,
\partial/\partial y_1, \ldots ,\partial/\partial y_8$ are
obviously independent so let $T$ be the orientation on $S^8\times
S^8$ determined by this local frame-field. The orientation $t$ on
$V = U_2^{\oplus 5}$ is determined by $15$ coordinate functions
$\{b_i^{++}, b_i^{+-}, b_i^{-+}\}_{i=1}^5$ where, as in
Section~\ref{sec:revisited}, $b_i^{\epsilon_1 \epsilon_2} =
\mu_i(H^{\epsilon_1, \epsilon_2})$ and $H^{\epsilon_1,
\epsilon_2}$ is the quadrant in $\mathbb{R}^8$ determined by the
oriented hyperplanes $H_1$ and $H_2$ corresponding to the signs
$\epsilon_1, \epsilon_2$.

\begin{figure}[hbt]
\centering
\includegraphics[scale=0.70]{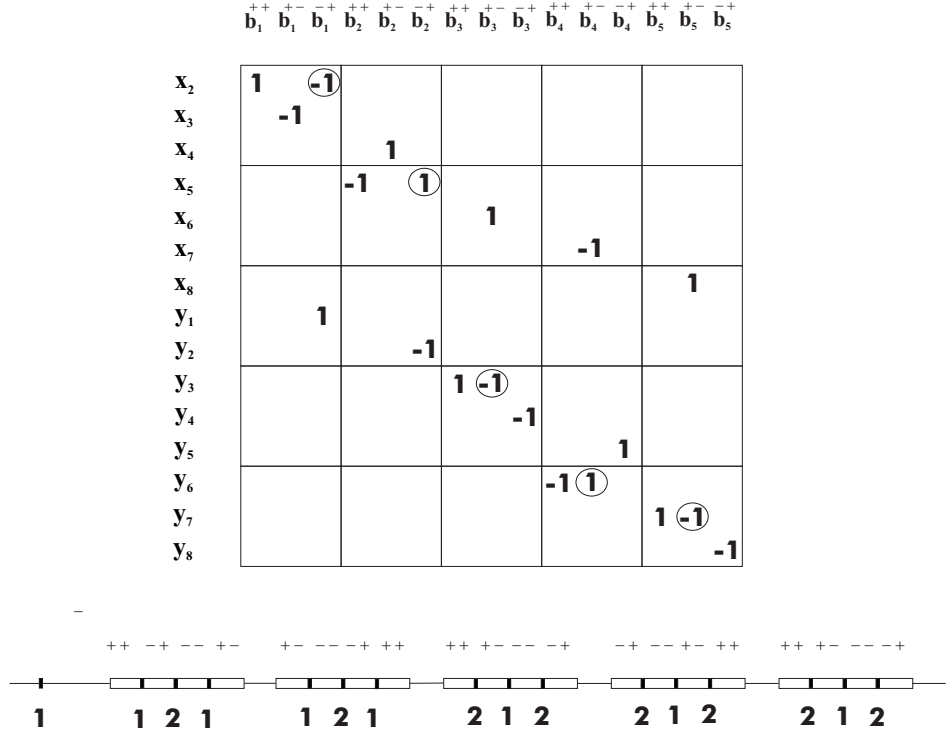} \caption{Jacobian matrix: $AAABBB$-case}
\label{fig:matrica1}
\end{figure}

\begin{figure}[hbt]
\centering
\includegraphics[scale=0.70]{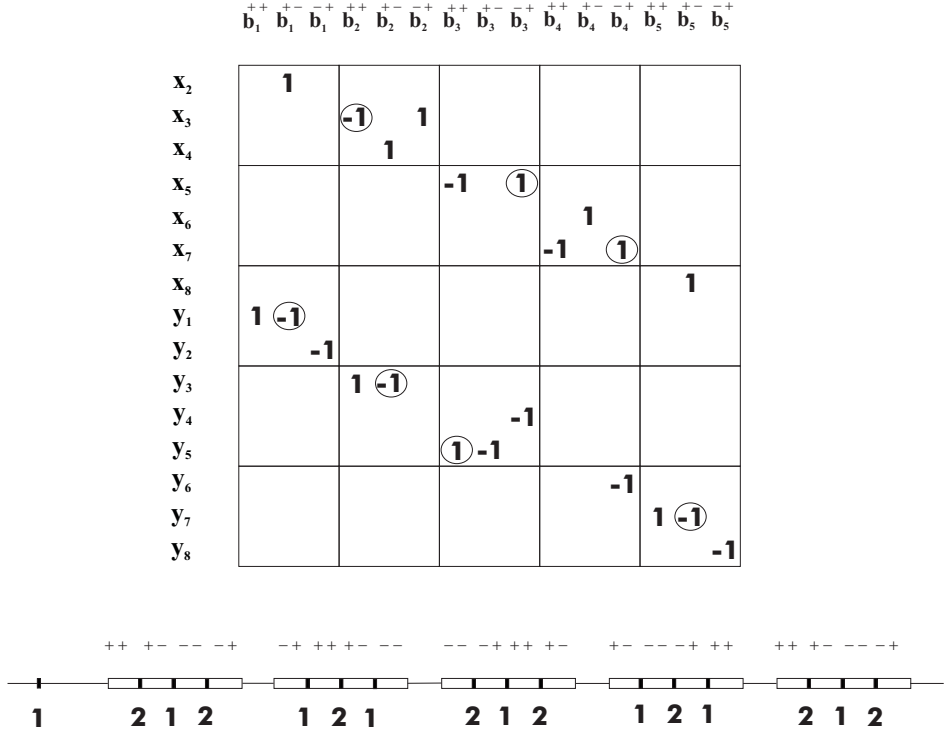} \caption{Jacobian matrix: $ABABAB$-case}
\label{fig:matrica2}
\end{figure}

Note that by fixing the coordinate $x_1$, the signed words
$AAABBB++$ and $ABABAB++$ can be interpreted as points on circles
$C_1$ and $C_2$ respectively. In both cases $x_2,\ldots, x_8,
y_1,\ldots , y_8$ can be chosen as the coordinates in the normal
slices to circles $C_1$ and $C_2$ at these points. Following the
procedure from Section~\ref{sec:algorithm} we should compare the
orientations in the normal bundles to these points with the
orientation $t$. This amounts to computing the signs of the
corresponding Jacobian matrices, Figures~\ref{fig:matrica1} and
\ref{fig:matrica2}. The circled entries in these matrices can be
removed by elementary column/row operations, so they have the same
determinant as the corresponding signed, permutation matrices.

It turns out that in the $AAABBB$ case the signed permutation has
precisely $7$ negative signs and $19$ transpositions while in the
$ABABAB$ case the corresponding sign permutation has $7$ negative
signs and $30$ transpositions. This implies that $C_1$ and $C_2$
have the opposite orientations in their normal bundles, hence the
elements $[M_1]$ and $[M_3]$ indeed cancel out. This completes the
proof of the theorem. \hfill $\square$

\subsubsection{The case $(6m+2, 4m+1, 2)$}

Recall (Definition~\ref{def:orbit}) that a generating class is a
$\mathbb{D}_8$-orbit of cyclic signed words while a cyclic signed
word is (Definition~\ref{sign-words}) an equivalence class of
signed words. We already know (Observation~\ref{observ1}) that
there is a close connection between cyclic signed words/generating
classes on one side and individual circles/minimal
$\mathbb{D}_8$-invariant submanifolds of the solution manifold on
the other. In particular, each generating class is a union of $2$,
$4$ or $8$ distinct cyclic signed words (its ``connected
components'').

The following definition introduces the key combinatorial
functions which permit us to translate the analysis given in
Section~\ref{sec:algorithm} into an algorithm for computing the
relevant obstruction class.

\begin{defin}
\label{def:funkc}  Let $\beta(n)$ be the number of generating
classes of cyclic signed words of length $2n$ which consist of $4$
different cyclic signed words. Suppose that $G$ is a generating
class of cyclic signed words consisting of precisely $2$ cyclic
signed words and let $\omega = C_1C_2\ldots C_{2n}++$ be the
lexicographically first signed $\{A,B\}$-word in this class.
Assume that this cyclic signed word is oriented according to the
action of the cyclic group $\mathbb{Z}/(2n)$ permuting its
letters. Let $\epsilon(G)$ be $+1$ or $-1$ depending on whether
the action of the generator $\gamma\alpha$, a ``rotation'' through
$90^\circ$, agrees with this orientation or not. Associate a
Jacobian matrix $J(\omega)$ to this word by a direct
generalization of the algorithm for constructing these matrices,
described in the proof of Theorem~\ref{thm:852} where it led to
Figures~\ref{fig:matrica1} and \ref{fig:matrica2}. Let $\eta(G)$
be the sign of the determinant of this matrix. Define $\alpha(n)$
(respectively $\beta(n)$) as the number of generating classes $G$
with two signed cyclic word components such that $\epsilon(G) =
\eta(G)$ (respectively $\epsilon(G) \neq \eta(G)$).
\end{defin}

\begin{rem}
Note that $\alpha(n), \beta(n), \gamma(n)$ are combinatorial
functions which have much in common with very well known functions
enumerating the number of cyclic words in a given alphabet.
However, an explicit formula for these functions, or at least for
$\alpha(n)$ and $\gamma(n)$ remains an interesting open problem.
\end{rem}

\medskip\noindent
{\bf Proof and comments on Theorem~\ref{thm:glavna}:} The proof of
Theorem~\ref{thm:glavna} follows step by step the procedure
outlined in the proof of the special case $(8,5,2)$. The
computation of the associated Jacobian matrices can be simplified
as follows. We illustrate the idea on the special case of the
matrix associated to the $AAABBB$-case of the solution manifold
associated to the triple $(8,5,2)$, Figure~\ref{fig:matrica1}.
Instead of working with $x_2,\ldots , x_8, y_1, \ldots, y_8$, as
local coordinates, we could put these functions in the order of
appearance, relative the equipartition of intervals shown in
Figure~\ref{fig:matrica1}. By inspection of
Figure~\ref{fig:matrica1} we observe that the natural order is
\[
x_2   \quad y_1 \quad x_3 \quad x_4 \quad   y_2 \quad x_5 \quad
y_3 \quad  x_6 \quad y_4 \quad y_5 \quad  x_7 \quad y_6 \quad y_7
\quad x_8 \quad y_8 .
\]
This system of functions has an advantage that the Jacobian matrix
with respect to this system is a $3\times 3$-block diagonal
matrix. Note that the sign of this matrix is equal to the sign of
the matrix on Figure~\ref{fig:matrica1}, multiplied by the sign of
the corresponding shuffle permutation, associated to the word
$AAABBB$. \hfill $\square$

\subsubsection{The case $(6m-1, 4m-1, 2)$}
\label{sec:general2}

For completeness we discuss here the case of a triples $(6m-1,
4m-1, 2)$ although these results have an alternative proof based
on completely different ideas, Section~\ref{sec:ideal}.

\begin{prop}
\label{prop:recporec} Suppose that $(d,j,k) = (6m-1,4m-1,2)$ where
$m$ is a positive integer. Then there exists an equipartition of
$j=4m-1$ mass distributions in $R^d = R^{6m-1}$ if an element $o =
o_1 + o_2 \in {\mathbb{Z}/2}\oplus {\mathbb{Z}/2}$ is nonzero,
where $o_1$ and $o_2$ are determined by the following congruences
modulo $2$,

\begin{equation}
\label{eqn:ooo} o_1 \equiv_2 O_1(m) := \sum_{\mbox{\scriptsize
$k\vert 2m$} \atop{\mbox{\scriptsize $k$ {\rm is odd} }}} A(k)
\qquad o_2 \equiv_2 O_2(m) := \sum_{\mbox {\scriptsize $k\vert
2m$} \atop{\mbox{\scriptsize
 $k$ {\rm is even} }}} A(k)
\end{equation}
where $A(k)$ is the number of $\ast$-primitive, circular words of
length $2k$, Definition~\ref{special}.
\end{prop}

\medskip\noindent
{\bf Proof:} The proof is a generalization of the analysis given
in Example~\ref{exam:circles}. In the case $(d,j,k) =
(6m-1,4m-1,2)$ the relevant obstruction lives in the group
$\mathbb{Z}/2\oplus\mathbb{Z}/2$. In order to compute this
obstruction we ``list'' all generating classes of cyclic signed
words of length $4m$ and determine the corresponding stabilizers,
which allows us to apply Proposition~\ref{recognition1}. Observe
that neither $(\alpha)$ nor $(\beta)$ appear as stabilizers of
connected components of generating classes of cyclic signed words.
On the other hand the groups $(\gamma), (\alpha\beta\gamma)$,
respectively $(\alpha\gamma), (\beta\gamma)$ do appear as
stabilizers in the generating classes corresponding to
$\ast$-primitive, circular words and it is not difficult to
distinguish these two cases. Suppose that $[w]$ is the generating
class of a special word $w = (aa^\ast)\ldots (aa^\ast)$ where
$aa^\ast$ is a $\ast$-primitive word of length $2k$. Then a
component of $[w]$ is stabilized by $(\gamma)$, respectively by
$(\alpha\gamma)$, depending on whether $k$ is even or odd.  In
light of Proposition~\ref{recognition1} this immediately leads to
formula (\ref{eqn:ooo}). \hfill $\square$.

\begin{cor}
\label{cor:general2}
\[
\Delta(2^{q+1}-1,2) = 3\cdot 2^q -1 .
\]
\end{cor}

\medskip\noindent
{\bf Proof:} It follows from Propositions~\ref{prop:recporec} and
\ref{prop:racun} that only in the case $m = 2^p$ the obstruction
element $o = o_1 + o_2 \in {\mathbb{Z}/2}\oplus {\mathbb{Z}/2}$ is
nonzero. The details are left to the reader. \hfill $\square$

\section{Cohomological Methods}\label{sec:ideal}

\subsection{Ideal valued cohomological index theory}
\label{cohomological index}

A standard tool for proving non existence of equivariant maps is a
cohomological index theory,  \cite{Do} \cite{Fa-Hu} \cite{I-R}
\cite{Ziv2}. A particularly useful form of this theory is the so
called ideal valued cohomological index theory developed by
E.Fadell and S.~Husseini \cite{FH2}, see also \cite{Jaw} and
\cite{Ziv3}.

In this section we demonstrate how this theory can be applied to
the equipartition problem and compare it with the obstruction
theory approach from previous sections.

\begin{theo} \label{thm:dickson} Let
\begin{equation}
\label{eq:dickson} {\mathbb P}_k = {\rm Det}\left[
\begin{array}{ccccc}
x_1 & x_1^2 & x_1^4 & \dots & x_1^{2^{k-1}}\\

x_2 & x_2^2 & x_2^4 & \dots & x_2^{2^{k-1}}\\

\vdots & \vdots & \vdots & \ddots & \vdots \\

x_k & x_k^2 & x_k^4 & \dots & x_k^{2^{k-1}}
\end{array}\right] \in\mathbb{F}_2[x_1,\ldots ,x_k]
\end{equation}
be a Dickson polynomial. Then $(d,j,k)$ is an {\em admissible}
triple if
\begin{equation}
\label{eq:FH} ({\mathbb P}_k)^j\notin{\rm
Ideal}\{x_1^{d+1},\ldots,x_k^{d+1}\} .
\end{equation}
\end{theo}

\medskip\noindent
{\bf Proof:} The proof is based on the ideal valued, cohomological
index theory, as developed by E.~Fadell and S.~Husseini,
\cite{FH2}. Recall that the equipartition problem can be reduced,
Section~\ref{sec:revisited}, to the question of the (non)existence
of a $W_k$-equivariant map $A : (S^d)^k\rightarrow S(U_k^{\oplus
j})$ where  $W_k = (\mathbb{Z}/2)^{\oplus k}\rtimes S_k$, $U_k$ is
a $W_k$-representation described in Section~\ref{sec:revisited}
and $S(V)$ is a unit sphere in $V$. The representation $U_k$ was
characterized by the property that its restriction on the subgroup
$H = (\mathbb{Z}/2)^{\oplus k}$ is equivalent to the regular
representation of the group $(\mathbb{Z}/2)^{\oplus k}$.

The cohomology of the classifying space $BH$ with
$\mathbb{F}_2$-coefficients is $H^{\ast}(BH;\mathbb{F}_2)=
\mathbb{F}_2 [x_1,...,x_k]$ and the corresponding indices are
$${\rm Index}^H\left(\left(S^d\right)^k\right)=\left(x_1^{d+1},...,x_k^{d+1}\right),$$

$${\rm Index}^H\left(S\left(U_k^{\oplus
j}\right)\right)=\left(\left(\mathbb{P}_k\left(x_1,...,x_k
\right)\right)^j\right),$$

\noindent where
 \[
 \mathbb{P}_k\left(x_1,...,x_k\right)=x_1\cdots
x_k\left(x_1+x_2\right) \cdots \left(x_{k-1}+x_k\right)\cdots
\left(x_1+\cdots +x_k \right)
 \] is a Dickson polynomial. It is well
known  \cite{Will} that $\mathbb{P}_k$ can be expressed in the
form of determinant (\ref{eq:dickson}), or more explicitly
 $$
 \mathbb{P}_k\left(x_1,...,x_k\right)=\sum_{\sigma \in \Sigma_k} x_{\sigma
(1)}^{2^{k-1}} x_{\sigma (2)}^{2^{k-2}}\cdots x_{\sigma (k)}.
 $$
 A sufficient condition, \cite{FH2}, for a nonexistence of a
 $H$-equivariant map $f : X\rightarrow Y$ is the relation
 ${\rm Index}^H(Y) \not\subseteq {\rm Index}^H(X)$. In our case
 this relation takes the form of the condition (\ref{eq:FH}) and
 the result follows. \hfill $\square$

\medskip

In  order to apply Theorem~\ref{thm:dickson} we search in
$\left(\mathbb{P}_k\left(x_1,...,x_k\right)\right)^j$ for a
summand where the biggest exponent is as small as possible. From
the properties of the binomial coefficients over $\mathbb{F}_2$ we
deduce that the summands we are looking for are of the form
$$ m=\left(x_1^{2^{k-1}}x_2^{2^{k-2}}\cdots x_k\right)^{2^q}\cdot
\left(x_1x_2^2\cdots x_k^{2^{k-1}}\right)^r,$$

\noindent where $j=2^q+r$ and $0\leq r\leq 2^q-1$. It is not
difficult to see that there exists a summand with the exponent
$2^{k+q-1}+r$. If $2^{k+q-1}+r\leq d$ then the condition
(\ref{eq:FH}) from Theorem~\ref{thm:dickson} is fulfilled and as a
consequence we obtain the following result.

\begin{theo}
\label{thm:index}
 \begin{equation}
\label{eq:bound}
  \Delta (2^q+r,k)\leq 2^{k+q-1}+r.
 \end{equation}
\end{theo}
It is interesting to compare the inequality (\ref{eq:bound}) with
the only existing general upper bound $\Delta (j,k)\leq j2^{k-1}$
from \cite{Ram}. Since
\[
\left(2^q+r\right)2^{k-1}-2^{k+q-1}-r=r\left(2^{k-1}-1\right)\geq
0,\]
our estimate is the same as the one from \cite{Ram} in the
case $r=0$ and strictly better in all other cases. In the special
case $k=2$ we obtain the inequalities
$$3\cdot 2^{q-1}+\frac {3r}2\leq \Delta (2^q+r,2)\leq 2^{q+1}+r.$$
Note that the best estimate is obtained when $j$ is slightly less
then a power of $2$. For example if $r=2^q-1,$ i.e. if
$j=2^{q+1}-1$, we obtain the exact value
$$\Delta (2^{q+1}-1,2)=3\cdot 2^q-1,$$
which is a result already obtained in Section~\ref{sec:general2}
(Corollary~\ref{cor:general2}). Also, we obtain almost precise
values for $\Delta (j,2)$ when $j=2^{q+1}-2$ or $j=2^{q+1}-3$

\begin{equation}
\label{eqn:nejednacine} 3\cdot 2^q-3\leq \Delta (2^{q+1}-2,2)\leq
3\cdot 2^q-2 \qquad \qquad 3\cdot 2^q-4\leq \Delta
(2^{q+1}-3,2)\leq 3\cdot 2^q-3.
\end{equation}
Note that we already know that $\Delta(2^{q+1}-2, 2) = 3\cdot 2^q
-3$ (Proposition~\ref{prop:case0}). Also note that the result
$\Delta(5,2) = 8$ (Theorem~\ref{thm:852}) cannot be obtained by
this method since the general upper bound (\ref{eqn:nejednacine})
implies only the inequality $\Delta(5,2)\leq 9$, already known to
Ramos \cite{Ram}.

Here are some particular examples which illustrate the power of
Theorem~\ref{thm:index}. The best previously known upper bounds
are given in parentheses.

\begin{equation}
\begin{array}{ccc}
17\leq \Delta (7,3)\leq 19  \; \; (28)  &&  14\leq \Delta
(6,3)\leq 18 \; \; (24)\\
35\leq \Delta (15,3)\leq 39 \; \; (60) && 33\leq \Delta
(14,3)\leq 38 \; \; (56) \\
27\leq \Delta (7,4)\leq 35 \; \; (56) && 23\leq \Delta (6,4)\leq
34 \; \; (48)\\
57\leq \Delta (15,4)\leq 71 \; \; (120) && 53 \leq \Delta
(14,4)\leq 70 \; \; (112)
\end{array}
\end{equation}

\end{document}